\documentclass[12pt,a4paper]{amsart}
\usepackage[english]{babel}
\usepackage[utf8]{inputenc}
\usepackage{comment}
\usepackage{amssymb}
\usepackage{amsmath} 
\usepackage{amsthm} 
\usepackage{enumitem}
\usepackage{bbm}
\usepackage{bbold}
\usepackage{todonotes}
\usepackage{verbatim}

\usepackage{xcolor}
\pagecolor{white}

\newtheorem{thm}{Theorem}[section]
\newtheorem{re}[thm]{Remark}

\newtheorem{lemma}[thm]{Lemma}

\numberwithin{equation}{section}

\DeclareMathOperator*{\res}{ res}

\DeclareMathOperator{\ord}{ ord}

\title{Averages of long Dirichlet polynomials with modular coefficients}
\author{Brian Conrey, Alessandro Fazzari}
\address{American Institute of Mathematics, 600 East Brokaw Road, San Jose, CA 95112, US}
\email{conrey@aimath.org, fazzari@aimath.org}

\begin{document}
\maketitle
\begin{abstract}
We study the moments of $L$-functions associated with primitive cusp forms, in the weight aspect. In particular, we obtain an asymptotic formula for the twisted moments of a \textit{long} Dirichlet polynomial with modular coefficients. This result, which is conditional on the Generalized Lindel\"of Hypothesis, agrees with the prediction of the recipe by Conrey, Farmer, Keating, Rubinstein and Snaith. 
\end{abstract}

\section{Introduction and statement of the main result}

Holomorphic modular forms, originally introduced in number theory for their connection to class field theory and elliptic curves, have nowadays applications in many directions. In particular, the central value of $L$-functions associated to modular forms has been studied extensively (see e.g. \cite{Conrey-Iwaniec, Iwaniec-Sarnak}), in view of its connection with deep arithmetic problems, such as the Birch and Swinnerton-Dyer conjecture. \par

In this paper, we focus on averages of Hecke $L$-functions at the central point\footnote{We normalize the $L$-functions so that the functional equation links values at $s$ to $1-s$, so the central point is $s=\frac{1}{2}$.}, in the weight aspect; in particular, we consider 
the normalized Hecke basis $H_k=H_k(1)$ of the space of holomorphic cusp forms of level 1 and (even) weight $k\geq 6$ 
and we study the $r$-th moment of the associated $L$-functions, i.e.
$$ M_r(k)=\sum_{f\in H_k}\!\! {\vphantom{\sum}}^{h}L_f(\tfrac{1}{2})^r,  $$
where the superscript $^h$ indicates the usual harmonic weight arising from Petersson norm and involving the symmetric square of the $L$-function. 
For small integer values of $r$, this problem has been studied in various works \cite{Balkanova-Frolenkov,Bykovskii-Frolenkov, Bykovskii-Frolenkov2, Balkanova-Frolenkov2, Balkanova-Frolenkov3, Frolenkov, Khan, HK} but for an arbitrary integer $r$ an asymptotic formula for $M_r(k)$ is currently out of reach. However, very precise conjectures have been formulated; $L$-functions associated to primitive cusp forms of weight $k$ form an orthogonal family in the sense of Katz and Sarnak \cite{Katz-Sarnak}, thus one expects $M_r(k)$ to be asymptotic, as $k$ goes to infinity, to a certain (explicit) polynomial in $\log k$ of degree $\frac{r(r-1)}{2}$. This conjecture comes from the now classical analogy between random matrix theory and number theory \cite{Montgomery, Katz-Sarnak, KS2000, Odlyzko}, which for many purposes allows us to predict the behavior of $L$-functions, by studying the classical group of matrices that models the family of $L$-functions we are working with, being typically unitary, symplectic, or orthogonal. In particular, by using a heuristic number theoretical machinery, Conrey, Farmer, Keating, Rubinstein and Snaith \cite{CFKRS} developed a recipe that conjecturally produces all the main and lower order terms for shifted moments of $L$-functions. Also in the case of modular forms, the recipe leads to an explicit conjecture for integral moments (see \cite[Conjecture 1.5.5, Conjecture 4.5.1]{CFKRS}). A compact way to write the expected result from the recipe is the following; 
given $A=\{\alpha_1,\dots,\alpha_r\}$ a set of small complex shifts, \cite{CFKRS} conjectures that
\begin{equation}\label{7june.1}
\sum_{f\in H_k}\!\! {\vphantom{\sum}}^{h} \prod_{\alpha\in A}L_f(\tfrac{1}{2}+\alpha)
\sim \sum_{V\subset A}i^{|V|k}
\left(\frac{k}{4\pi}\right)^{-2\sum_{\alpha\in V}\alpha}G_{A\smallsetminus V\cup V^-}(\tfrac{1}{2};1)
\end{equation}
where $V^-:=\{-v:v\in V\}$ and 
\begin{equation}\label{G_l(A)}
G_A(s;l) := \sum_{m=1}^{\infty} \frac{1}{m^s}\sum_{m_1\cdots m_r=m}\frac{\prod_p\frac{2}{\pi}\int_0^{\pi} \big(\prod_{i=1}^r U_{\ord_p(m_i)}(\cos\theta)\big) U_{\ord_p(l)}(\cos\theta)\sin^2\theta\;d\theta}{m_1^{\alpha_1}\cdots m_r^{\alpha_r}}
\end{equation}
with $U_m$ denoting the Chebyshev polynomials, defined in \eqref{defCP}. 
See \cite{CFKRS} for the analytic continuation properties of $G_{A}(s;l)$ at $s=\frac{1}{2}$.
Equation \eqref{7june.1} shows the structure of the shifted $r$-moment of $L_f(\tfrac{1}{2})$; each term of the outer sum over the subsets of $A$ gives a main term, which is called an \lq\lq$\ell$-swap\rq\rq  \;term when $V$ has cardinality $\ell$. Moreover, the so-called 0-swap term, i.e. the one corresponding to $V=\emptyset$, 
lays the groundwork for the rest of the formula,
as all the other terms can be essentially obtained by changing sign of some of the shifts.\par

A classical approach to attack \eqref{7june.1} is to truncate the Dirichlet series of the product of $L$-functions on the left hand-side. Therefore, studying averages of Dirichlet polynomials that approximate $L$-functions is a central problem in the context of moments. By the Conrey-Keating heuristic \cite{CK1, CK2, CK3, CK4, CK5}, we expect the length $X$ of the polynomial approximation to play a crucial role; in particular, if we approximate the Dirichlet series of $\prod_{\alpha\in A}L_f(\tfrac{1}{2}+\alpha)$ with a polynomial of length $k^{2L}<X<k^{2(L+1)}$ (note that $k^2$ is the size of the conductor of $L_f$), then only the $\ell$-swap terms with $\ell\leq L$ from the recipe actually contribute to the main term (see also \cite{Bettin-Conrey, Sieg-Caroline} for more details on how the length of the polynomial approximation influences the final result). For example, if $X<k^2$ then only the 0-swap term gives a non-negligible contribution, while all the others contribute $o(1)$. In this paper, we treat the so-called \textit{long Dirichlet polynomial approximation}, i.e. the case $k^2<X<k^4$, for which we expect that only the 0-swap and 1-swap terms contribute to the main term and the remaining terms are negligible.

However, the current technology in the literature does not allow to go so far with moments of $L$-functions associated with cusp forms. For example, an asymptotic formula has been proven only in the case of the first and second moments, but already the third moment is not known. In the language we used above, we cannot go beyond length $Xl\ll k^2$ and in this range the 1-swaps do not appear yet. As we will describe in detail below, the difficulty is  given by the $J$-Bessel function coming from the off-diagonal term in Petersson formula. 
The nature of this limitation is the same that appears in \cite[Theorem 1.2]{ILS}, where Iwaniec, Luo and Sarnak compute the one-level density for complex zeros of modular $L$-functions only for test function with Fourier support in $(-1,1)$ instead of the \lq\lq typical\rq\rq $(-2,2)$ \cite[Theorem 1.2]{ILS}. In the same paper, they overcome this issue\footnote{As Iwaniec, Luo and Sarnak mention below Thereom 1.2 in \cite{ILS}, this issue is due to the fact that the conductor of the $L$-function in this family $k^2$ is twice  as large as the cardinality of the family itself, on a logarithmic scale.} at the cost of performing an extra average over $k$, that allows them to enlarge the admissible Fourier support of the test function from $(-1,1)$ to $(-2,2)$. Here, we will do the same and, by adding an extra average over the weight $k\asymp K$, we will \lq\lq double\rq\rq the admissible length of the polynomial to $Xl\ll K^4$.
In terms of the Bessel function $J_{k-1}(x)$, that behaves roughly like $e^{ix}/\sqrt{x}$ for $x\gg k$, the advantage is that its oscillating nature can be seen much more precisely once we average over $k$ (see \eqref{20mar.3} for the formal statement and the beginning of Section 8 in \cite{ILS} for a comprehensive account).

The following is the main result of this paper, proven conditionally on the assumption of the Generalized Lindel\"of Hypothesis for Dirichlet $L$-functions (GLH).

\begin{thm}\label{thm}
Assume GLH and let $K\geq 1$. For any given positive integer $r$, let $A=\{\alpha_1,\dots,\alpha_r\}$ be a set of shifts of size $\ll (\log K)^{-1}$, 
$X>0$ and $l\in \mathbb N$ two parameters and $\psi$ a smooth function supported in $(0,1)$.
Let also $h$ be a smooth function, compactly supported on the positive real numbers. Then, as $K\to\infty$, we have
\begin{equation}\begin{split}\label{24mar.1}
&\frac{4}{HK}\sum_{k\equiv0\;(4)}h\left(\frac{k-1}{K}\right)\sum_{f\in H_k}\!\! {\vphantom{\sum}}^{h} \lambda_f(l)\sum_{n=1}^\infty \psi\left(\frac{n}{X}\right)\frac{1}{\sqrt{n}}\sum_{n_1\cdots n_r=n}\frac{\lambda_f(n_1)\cdots \lambda_f(n_r)}{n_1^{\alpha_1}\cdots n_r^{\alpha_r}} \\
&= 
\frac{4}{HK}\sum_{k\equiv0\;(4)}h\left(\frac{k-1}{K}\right)
\frac{1}{2\pi i}\int_{(\varepsilon)}
%\frac{X^z}{z}
\tilde\psi(z)X^z
\sum_{\substack{V\subset A \\ 0\leq |V|\leq 1}}
 \left(\frac{k}{4\pi}\right)^{-2\sum_{\alpha\in V}(\alpha+z)}G_{A_z\smallsetminus V_z\cup V_z^-}(\tfrac{1}{2};l)dz \\
&\hspace{11.3cm}+ O_{\varepsilon}\bigg((XlK)^{3r\varepsilon}\frac{\sqrt{Xl}}{K^2}\bigg),
\end{split}\end{equation}
for every $\varepsilon>0$, with $V^-:=\{-v:v\in V\}$, $V_z:=\{v+z:v\in V\}$ and $G_l(A)$ as in \eqref{G_l(A)}.
\end{thm}

We notice that we average over $k\equiv 0\;(4)$, to avoid the trivial cancellation of the odd (i.e. $|V|$ even) swaps, given by the factor $i^{|V|k}$ in \eqref{7june.1} that comes from the $(-1)^{k/2}$ factor in the functional equation of modular $L$-functions \eqref{feL_f}. Also, on the right hand side of \eqref{24mar.1} the sum over $k$ can be explicitly computed, since $G_A$ does not depend on $k$.

Note that, in the range $Xl\asymp K^\eta$ with $2<\eta<4$, the error term above is $o(1)$ (if $\varepsilon$ is sufficiently small), while the main terms are of size $\gg 1$.
The above theorem then shows that, under GLH, the 0-swap and 1-swap terms predicted by the recipe are correct for the orthogonal family of $L$-functions associated to primitive cusp forms of weight $k$, when we average over $k\asymp K$. Results of the same nature have been proven also for other families, in particular by Baluyot and Turnage-Butterbaugh \cite{Sieg-Caroline} for the unitary family of Dirichlet $L$-functions of primitive characters and by Conrey and Rodgers \cite{Brad-Brian} for the symplectic case of quadratic Dirichlet $L$-functions; both these results have been proven under GLH. We also mention \cite{Hamieh-Ng}, where Hamieh and Ng study the (continuous) family of $\zeta$ in the $t$-aspect, assuming the expected asymptotic formula for correlations of divisor sums. A motivation for all these results lies in the Conrey-Keating heuristic \cite{CK5}, which suggests how to ``patch together'' the 1-swap terms in order to construct the general $\ell$-swap term (see \cite{Brian-Sieg} for a complete overview on the topic). This approach actually requires a twisted version of the moments of the long Dirichlet approximation, which has therefore been included in Theorem \ref{thm}. \par

The proof of Theorem \ref{thm} can be summarized as follows. Section \ref{S2} is devoted to some arithmetical manipulations, which allow us to write in a compact way an arbitrary convolution of modular coefficients $\lambda$. Then, in Section \ref{S3}, we apply our main device, being the Petersson trace formula \eqref{PeterssonFormula}; the ``$\delta$-term'' from the trace formula can be easily matched with the 0-swap term from the recipe. The ``off-diagonal term'' brings into play Kloosterman sums, which lead us to study a generalized Estermann function, i.e. an exponential twist of the initial $L$-function. After studying the polar structure of the Estermann function (see Lemma \ref{lemmaEF} and Remark \ref{RemarkEstermann}), by a contour shift, we extract the contributions from the poles and bound the remaining integral, relying on GLH; this is achieved in Section \ref{S5}. Finally %, with some gymnastic, 
Section \ref{S6} proves that the residues of the Estermann function give precisely the expected 1-swap terms from the recipe, concluding the proof.\par

%For longer polynomials, i.e. when $X\asymp K^\eta$ with $\eta>4$, to get an asymptotic formula for the $r$-th moment, one would need to detect also higher swaps. For instance, the problem of detecting the 2-swap terms boils down to a Voronoi summation formula for the Estermann function associated to the $r$-fold divisor function. \textcolor{red}{remove?}\par
%
%Finally, we note that Theorem \ref{thm} is known unconditionally for $r=1,2$. Moreover, without assuming GLH, one can prove the same result but for shorter polynomials. More precisely, Theorem \ref{thm} holds unconditionally for $\eta < \frac{4+2r\theta_B}{1+r\theta_B}$, where $\theta_B=\frac{3}{16}$ denotes the exponent in the Burgess' bound \cite{Burgess}, $L(\frac{1}{2},\chi)\ll q^{\theta_B+\varepsilon}$, for $\chi$ character (mod $q$). \textcolor{red}{remove?}

\section{Arbitrary convolutions of modular coefficients}\label{S2}

Let $k\geq 6$ be an even positive integer and $H_k=H_k(1)$ be the set of primitive cusp forms of weight $k$ and level 1, eigenfunctions  of the Hecke operators. 
The $L$-function associated with $f\in H_k$ is 
$$ L_f(s) = \sum_{n=1}^\infty\frac{\lambda(n)}{n^s} = \prod_p\bigg( 1-\frac{\lambda(p)}{p^s} + \frac{1}{p^{2s}} \bigg)^{-1}, \quad \Re(s)>1$$
where $\lambda(n) = \lambda_f(n)$ is the (normalized) $n$-th Fourier coefficient of the cusp form $f$, given by the expansion
$$ f(z)=\sum_{n=1}^\infty \lambda(n)n^{(k-1)/2}e^{2\pi i n z} .$$
This $L$-function satisfies the functional equation 
\begin{equation}\label{feL_f} 
\Lambda_f(s)=(2\pi )^{-s}\Gamma\left(s+\frac{k-1}{2}\right)L_f(s)= (-1)^{k/2}\Lambda_f(1-s). 
\end{equation}
See \cite{IwaniecTopics} for a more comprehensive account about modular forms.

When we look at the $r$-th moment, denoting by $A=\{\alpha_1,\dots,\alpha_r\}$ a set of small shifts of size $\ll (\log k)^{-1}$, we need to write in a convenient way a convolution of modular coefficients. Namely, we have
\begin{equation}\begin{split}\label{5apr.0} 
%L_f(s+\alpha_1) \cdots L_f(s+\alpha_r) 
\prod_{\alpha\in A}L_f(s+\alpha)
= \sum_{m_1,\dots,m_r=1}^{\infty} \frac{\lambda(m_1)\cdots\lambda(m_r)}{m_1^{s+\alpha_1}\cdots m_r^{s+\alpha_r}} 
= \prod_p\sum_{m_1,\dots,m_r=0}^{\infty}\frac{\lambda(p^{m_1})\cdots\lambda(p^{m_r})}{p^{m_1(s+\alpha_1)+\dots+m_k(s+\alpha_r)}}
\end{split}\end{equation}
and we want to use the Hecke relation to write in a compact way the product of all the $\lambda$s at the numerator.
For instance, for $r=2$, the Hecke multiplicativity relation
\begin{equation}\label{HeckeRelation}\notag
\lambda(n)\lambda(m) = \sum_{d|(m,n)}\lambda\left(\frac{mn}{d^2}\right)
\end{equation}
allows us to write
\begin{equation}\label{Heckex2} 
L_f(s+\alpha)L_f(s+\beta) = \zeta(2s+\alpha+\beta) \sum_l\frac{\lambda(l)\tau_{\alpha,\beta}(l)}{l^s},
\end{equation}
where $\tau_{\alpha,\beta}(l) = \sum_{ab=l}a^{-\alpha}b^{-\beta}$ denotes the 2-fold shifted divisor function.
A similar formula can be proved also for $r=3$. Already for $r=4$, however, the arithmetic becomes quite hard; therefore, to deal with the general case $r\in\mathbb N$, we use a more succinct approach and we prove the following.

\begin{lemma}
In the above notations, for $\Re(s)>1$, we have
\begin{equation}\begin{split}\label{5apr.11} 
\prod_{\alpha\in A}L_f&(s+\alpha) 
= \sum_{m=1}^\infty\lambda(m)F_A(m,s),
\end{split}\end{equation} 
where
\begin{equation}\label{5apr.10} 
F_{A}(m,s) = \prod_{p^l||m}F_{A,p}(l,s)
\end{equation}
and
\begin{equation}\begin{split}\label{5apr.9}\notag
F_{A,p}(l,s) &= \sum_{m_1,\dots,m_r=0}^\infty \frac{c_l(m_1,\dots,m_r)}{p^{m_1(s+\alpha_1)+\dots+m_r(s+\alpha_r)}},\\
c_l(m_1,\dots,m_r)  &= \frac{2}{\pi}\int_0^\pi U_{m_1}(\cos\theta)\cdots U_{m_r}(\cos\theta) U_{l}(\cos\theta) \sin^2\theta\;d\theta.
\end{split}\end{equation}
\end{lemma}

Equation \eqref{5apr.11} provides us with a convenient expression for a product of $r$ copies of $L$-functions. In particular, we notice that the dependence on the cusp form $f$ is all encoded in the factor $\lambda(m)$, while the factor $F_A(m,s)$ is independent of $f$ and takes care of the arithmetic difficulties that rise when applying Hecke multiplicativity several times.

\proof
First, we notice that (by Euler product)
\begin{equation}\label{5apr.1}
\lambda(p^m) 
= \frac{1}{m!}\partial_z^m\left[\left(1-\lambda(p)z+z^2\right)^{-1}\right]_{z=0}.
\end{equation}
Moreover, by the Deligne bound $|\lambda(p)|\leq 2$, we can write
\begin{equation}\label{5apr.2}
\lambda(p) = 2\cos\theta_p. 
\end{equation}
For any nonnegative integer $m$, we introduce the Chebyshev polynomial $U_m$, defined by
\begin{equation}\label{defCP} U_m(\cos\theta)\sin\theta=\sin((m+1)\theta)\end{equation}
whose generating function is
\begin{equation}\label{5apr.3}
\sum_{j=0}^\infty U_j(t)X^j = \left( 1-2Xt+X^2 \right)^{-1}.
\end{equation}
In view of Equations \eqref{5apr.1}, \eqref{5apr.2} and \eqref{5apr.3}, we can identify
\begin{equation}\label{5apr.4} \lambda(p^m) = U_m(\cos\theta_p). \end{equation}
It is well-know that the Chebyshev polynomials form an orthonormal basis for the vector space of all the polynomials, with respect to the Sato-Tate measure, i.e. the inner product
$$ \frac{2}{\pi}\int_0^\pi g(\theta)h(\theta)\sin^2\theta\; d\theta. $$
Hence we have
\begin{equation}\label{5apr.5} 
U_{m_1}(\cos\theta)\cdots U_{m_r}(\cos\theta) = \sum_{l=0}^\infty c_l(m_1,\dots,m_r) U_l(\cos\theta)
\end{equation}
(note that the above is in fact a finite linear combination)
and the coefficients can be easily found by orthonormality, being
\begin{equation}\begin{split}\label{5apr.7} \notag
c_l(m_1,\dots,m_r)  = \frac{2}{\pi}\int_0^\pi &U_{m_1}(\cos\theta)\cdots U_{m_r}(\cos\theta) U_{l}(\cos\theta) \sin^2\theta\;d\theta.
\end{split}\end{equation}
Recalling \eqref{5apr.4}, Equation \eqref{5apr.5} gives a compact expression for the product of many $\lambda$ coefficients; more specifically, we have
\begin{equation}\label{5apr.6} \notag
\lambda(p^{m_1})\cdots \lambda(p^{m_r}) = \sum_{l=0}^\infty c_l(m_1,\dots,m_r) \lambda(p^{l}).
\end{equation}
Therefore, 
in view of Equation \eqref{5apr.0}, we have
\begin{equation}\label{5apr.8} \notag
\prod_{\alpha\in A}L_f(s+\alpha) = \prod_p \sum_{l=0}^\infty \lambda(p^l)\sum_{m_1,\dots,m_r=0}^\infty\frac{c_l(m_1,\dots,m_r)}{p^{m_1(s+\alpha_1)+\dots+m_r(s+\alpha_r)}},
\end{equation}
which proves the claim by definition of $F_A(m,s)$ and $F_{A,p}(l,s)$.
\endproof

Denoting
\begin{equation}\label{16june.1} C(X;\theta) = \left(1-2X\cos\theta+X^2\right)^{-1} \end{equation}
we now prove a lemma, which provides a further understanding of the behavior of the function $F_A(m,s)$, showing that in first approximation it behaves like
$$ \prod_{1\leq i<j\leq r}\zeta(2s+\alpha_i+\alpha_j) \times\frac{\tau_A(m)}{m^s}, $$
in analogy to \eqref{Heckex2}, where $\tau_A(m)$ denotes the $m$-th coefficient of the Dirichlet series of $ \prod_{\alpha\in A}\zeta(s+\alpha)$.
In addition, the following lemma studies the Dirichlet series associated with $F_A(m,s)$, twisted by any Dirichlet character. 

\begin{lemma}\label{lemmaEF}
For $F_A(m,s)$ and $C(X;\theta)$ defined in \eqref{5apr.10} and \eqref{16june.1} respectively, $\Re(s)>1$, $m\in\mathbb N$, we have
\begin{equation}\begin{split}\label{5apr.12}
F_A(m,s) 
= \prod_{p^\nu||m} \frac{2}{\pi}\int_0^\pi \prod_{\alpha\in A} C\left(\tfrac{1}{p^{s+\alpha}};\theta\right) \times U_\nu(\cos\theta) \sin^2\theta\;d\theta
\end{split}\end{equation}
and, for $p\nmid m$,
\begin{equation}\begin{split}\label{5apr.12added}
\frac{2}{\pi}\int_0^\pi& \prod_{\alpha\in A} C\left(\frac{1}{p^{s+\alpha}};\theta\right) \times U_{0}(\cos\theta)\sin^2\theta\;d\theta =1 + \sum_{1\leq i<j\leq r} \frac{1}{p^{2s+\alpha_i+\alpha_j}}+ O\left(\frac{1}{p^{3/2}}\right),
\end{split}\end{equation}
which provides an analytic continuation for $F_A(m,s)$ in the region $\Re(s)>\frac{1}{2}$.

Moreover, given $\chi$ a Dirichlet character, for $\Re(s)>1$ and $\Re(w)>1$, then
\begin{equation}\begin{split}\label{26june.1}
\mathcal F_A(w,s,\chi)
:&=\sum_{m=1}^{\infty} \frac{F_A(m,s)\chi(m)}{m^w} \\
&= \prod_p \frac{2}{\pi}\int_0^\pi \prod_{\alpha\in A} C\left(\frac{1}{p^{s+\alpha}};\theta\right) \times C\left(\frac{\chi(p)}{p^w};\theta\right)\sin^2\theta\;d\theta 
\end{split}\end{equation}
and
\begin{equation}\begin{split}\label{18july.1}
\frac{2}{\pi}\int_0^\pi &\prod_{\alpha\in A} C\left(\frac{1}{p^{s+\alpha}};\theta\right) \times C\left(\frac{\chi(p)}{p^w};\theta\right)\sin^2\theta\;d\theta \\
&=1 + \sum_{i=1}^{r}\frac{\chi(p)}{p^{s+w+\alpha_i}} + \sum_{1\leq i<j\leq r} \frac{1}{p^{2s+\alpha_i+\alpha_j}} + \sum_{1\leq i\leq j\leq r} \frac{\chi(p)^2}{p^{2s+2w+\alpha_i+\alpha_j}}  + O\left(\frac{1}{p^{3/2}}\right),
\end{split}\end{equation}
which gives a meromorphic continuation for $\mathcal F_A(w,s,\chi)$ in the region $\Re(s+w)>\frac{1}{2}$, $\Re(s)>\frac{1}{2}$.
\end{lemma}

\proof
We prove Equation \eqref{5apr.12} by direct computation; indeed for $\Re(s)>1$ 
\begin{equation}\begin{split}\notag
F_A(m,s) &= \prod_{p^\nu||m} F_{A,p}(\nu,s)
= \prod_{p^\nu||m} \sum_{m_1,\dots,m_r=0}^\infty \frac{c_\nu(m_1,\dots,m_r)}{p^{m_1(s+\alpha_1)+\dots+m_r(s+\alpha_r)}} \\
&= \prod_{p^\nu||m} \frac{2}{\pi}\int_0^\pi \prod_{j=1}^r\bigg(\sum_{m_j=0}^\infty U_{m_j}(\cos\theta)\frac{1}{p^{m_j(s+\alpha_j)}}\bigg)\times U_\nu(\cos\theta)\sin^2\theta\;d\theta  \\
&= \prod_{p^\nu||m} \frac{2}{\pi}\int_0^\pi \prod_{j=1}^r C\left(\tfrac{1}{p^{s+\alpha_j}};\theta\right)\times U_\nu(\cos\theta)\sin^2\theta\;d\theta,
\end{split}\end{equation}
by \eqref{5apr.3} and \eqref{16june.1}.
Equation \eqref{26june.1} can be proved similarly for $\Re(s)>1$ and $\Re(w)>1$, as
\begin{equation}\begin{split}\notag
\sum_{m=1}^{\infty} \frac{F_A(m,s)\chi(m)}{m^{w}} &
= \prod_p \sum_{l=0}^{\infty} \frac{F_{A,p}(l,s)\chi(p)^l}{p^{lw}}.
\end{split}\end{equation}
To put in evidence the polar structure of the above, we start by noticing that for $\Re(s)>0$ and $\Re(w)> 0$ (remember that the shifts $\alpha_1,\dots,\alpha_r$ are $\ll(\log k)^{-1}$)
\begin{equation}\label{july12.1}\notag
C\left(\frac{1}{p^{s+\alpha}};\theta\right) 
= 1+\frac{2\cos\theta}{p^{s+\alpha}}+\frac{4\cos^2\theta-1}{p^{2s+2\alpha}}+O\left(\frac{1}{p^{3/2}}\right) 
\end{equation}
and then
\begin{equation}\begin{split}\notag
\prod_{\alpha\in A} C\left(\frac{1}{p^{s+\alpha}};\theta\right) 
= 1 + \sum_{i=1}^{r}\frac{2\cos\theta}{p^{s+\alpha_i}} + \sum_{1\leq i<j\leq r} \frac{4\cos^2\theta}{p^{2s+\alpha_i+\alpha_j}} + \sum_{i=1}^r \frac{4\cos^2\theta-1}{p^{2s+2\alpha_i}} + O\left(\frac{1}{p^{3/2}}\right).
\end{split}\end{equation}
Moreover, for $|Y|<1$, we have
\begin{equation}\begin{split}\label{integraliutili}\notag
&\frac{2}{\pi}\int_0^\pi C(Y;\theta)\sin^2\theta\;d\theta = 1, 
\hspace{1.8cm} \frac{2}{\pi}\int_0^\pi 4\cos^2\theta \; C(Y;\theta)\sin^2\theta\;d\theta = 1+Y^2, \\
&\frac{2}{\pi}\int_0^\pi 2\cos\theta \; C(Y;\theta)\sin^2\theta\;d\theta = Y,
\hspace{0.5cm} \frac{2}{\pi}\int_0^\pi (4\cos^2\theta-1) \; C(Y;\theta)\sin^2\theta\;d\theta = Y^2.
\end{split}\end{equation}
Therefore, being $|C(Y;\theta)|\ll 1$ for $|Y|<1$, we get
\begin{equation}\begin{split}\notag
\frac{2}{\pi}\int_0^\pi& \prod_{j=1}^r C\left(\frac{1}{p^{s+\alpha_j}};\theta\right) C\left(\frac{\chi(p)}{p^{w}};\theta\right) \sin^2\theta\;d\theta 
\\& = 1 + \sum_{i=1}^{r}\frac{\tfrac{\chi(p)}{p^{w}}}{p^{s+\alpha_i}} + \sum_{1\leq i<j\leq r} \frac{1+(\tfrac{\chi(p)}{p^{w}})^2}{p^{2s+\alpha_i+\alpha_j}} + \sum_{i=1}^r \frac{(\tfrac{\chi(p)}{p^{w}})^2}{p^{2s+2\alpha_i}} + O\left(\frac{1}{p^{3/2}}\right) 
\end{split}\end{equation}
which concludes the proof of \eqref{18july.1}. Equation \eqref{5apr.12added} can be proved by a similar (and easier) computation.
\endproof

In the following result, we slightly refine the above lemma, including some arithmetical extra modifications that will appear in the following sections.

\begin{lemma}\label{lemmaEFrefined}
Let $c$ be an integer, $g$ a divisor of $c$, and $\chi$ a Dirichlet character modulo $q$. Then 
$$\mathcal F_A(w,s,\chi,c,g) :=\sum_{(m,\frac{c}{g})=1}\frac{F_A(mg,s)\chi(m)}{m^w}$$ 
has the same polar structure as 
$\mathcal F_A(w,s,\chi)$ from Lemma \ref{lemmaEF}. 
Moreover,
$\mathcal F_A(w,s;\chi,c,g)$ 
admits a meromorphic continuation in the region $\Re(s+w)>\tfrac{1}{2}$ and $\Re(s)>\frac{1}{2}$, where, under GLH for Dirichlet $L$-functions, we have
$$\mathcal F_A(w,s,\chi,c,g) \ll (qc)^{r\varepsilon}$$ for any $\varepsilon>0$.
\end{lemma}

\proof
Denoting $\nu_p(g)$ the $p$-adic evaluation of $g$, being $F_A(mg, s)/F_A(g, s)$ multiplicative in $m$, with some standard manipulations we can write
\begin{equation}\notag
\sum_{(m,\frac{c}{g})=1}\frac{F_A(mg,s)\chi(m)}{m^w}
= \sum_{m=1}^\infty \frac{F_A(m,s)\chi(m)}{m^w}\times \mathcal A_{c,g}(s,w)^{-1} \times \mathcal B_{c,g}(s) \times  \mathcal C_{c,g}(s,w),
\end{equation}
with
\begin{equation}\begin{split}\notag
\mathcal A_{c,g}(s,w) :& 
= \prod_{p| \frac{c}{g} }\sum_{m=0}^\infty\frac{F_{A,p}(m,s)\chi(p^m)}{p^{mw}},
\quad
\mathcal B_{c,g}(s) :=  \prod_{p|(\frac{c}{g},g)}F_{A,p}(\nu_p(g),s), \\
\mathcal C_{c,g}(s,w) :& =\prod_{\substack{p\nmid \frac{c}{g} \\ p| g}}\frac{\sum_{m=0}^\infty\chi(p^m) F_{A,p}(m+\nu_p(g),s)p^{-mw}}{\sum_{m=0}^\infty \chi(p^m)F_{A,p}(m,s)p^{-mw}} .
\end{split}\end{equation}
For $\Re(s)>\tfrac{1}{2}$ and $\Re(s+w)>\frac{1}{2}$, Lemma \ref{lemmaEF} (together with GLH) implies that $\mathcal F_A(w;s,\chi) \ll q^{r\varepsilon}$ for any $\varepsilon>0$. Moreover, in the same range for $s$ and $w$, we have
\begin{equation}\notag
\mathcal A_{c,g}(s,w)^{-1} \times \mathcal B_{c,g}(s) \times  \mathcal C_{c,g}(s,w) = \prod_{p|c} \bigg( 1+O\bigg(\frac{r}{\sqrt{p}}\bigg) \bigg).
\end{equation}
Then, since $1+O(\frac{r}{\sqrt{p}}) = e^{r/\sqrt{p}}$, the claim follows by noticing that, in the worst case scenario $c=\prod_{p\leq x}p$, we have
\begin{equation}\begin{split}\notag
\prod_{p|c}\exp\bigg(\frac{r}{\sqrt{p}}\bigg) 
&\leq \exp\bigg(\sum_{p\leq x} \frac{r}{\sqrt{p}}\bigg)  
\ll \exp\bigg(\frac{r\sqrt{\log c}}{\log\log c}\bigg) \ll c^{r\varepsilon}
\end{split}\end{equation}
by prime number theorem, being $\log c=\sum_{p\leq x}\log p \sim x$. 
\endproof

\section{Moments of a long Dirichlet polynomial approximation}\label{S3}

We denote by $\lambda_A(n)$ the  shifted convolution of copies of $\lambda$, namely 
$$ \prod_{\alpha\in A}L_f(s+\alpha) 
= \sum_{n_1,\dots,n_r=1}^{\infty}\frac{\lambda(n_1)\cdots\lambda(n_r)}{n_1^{s+\alpha_1}\dots n_r^{s+\alpha_r}} 
= \sum_{n=1}^{\infty}\frac{\lambda_A(n)}{n^s}. $$
Let's consider a positive integer $l$, a parameter $X>0$, a smooth function $\psi$ supported in $(0,1)$, a smooth function $h$ supported in $(h_1,h_2)$ (with $0<h_1<h_2<\infty$) and $\hat h(0)=H$, and $K$ a large parameter parameter. For simplicity we also assume $l,X\ll K^C$ for some fixed $C>0$, so that $l^\varepsilon$ and $X^{\varepsilon}$ are also $\ll K^\varepsilon$. Then we look at
\begin{equation}\label{20mar.6}\notag
M_A(l,X) := \frac{4}{HK}\sum_{k\equiv 0\;(4)}  h\left(\frac{k-1}{K}\right)\sum_{f\in H_k}\!\! {\vphantom{\sum}}^{h} \lambda(l) \sum_{n=1}^\infty \psi\left(\frac{n}{X}\right) \frac{\lambda_A(n)}{\sqrt{n}}, 
\end{equation}
where $\sum_{f\in H_k}\!\!\!\!\!\!\! {\vphantom{\sum}}^{h}\;\;*=\sum_{f\in H_k}*/\langle f,f\rangle$ is the usual average over the family of cusp forms, weighted by the Petersson norm, while the outer sum consists of an average over the weight $k\asymp K$.

\begin{re}\label{RemTotalMeasure}
Note that, by using Poisson summation formula and integrating by parts many times, one easily sees that 
$$ \frac{4}{HK}\sum_{k\equiv 0\;(4)}  h\left(\frac{k-1}{K}\right) = 1+O(K^{-A})$$ 
for any $A>0$. 
%\textcolor{blue}{
%\begin{equation}\begin{split}\notag
%\frac{4}{HK}\sum_{k\equiv 0\;(4)} h\left(\frac{k-1}{K}\right)
%&=\frac{4}{HK}\sum_m h\left(\frac{4m-1}{K}\right)
%=\frac{4}{HK}\sum_n \int_0^\infty h\left(\frac{4x-1}{K}\right)e^{-2\pi i xn}dx\\
%&=\frac{1}{H}\sum_n \int_{-\frac{1}{K}}^\infty h(y)e^{-\frac{\pi i}{4} (Ky+1)n}dy
%= 1+\frac{1}{H}\sum_{n\not=0} \int_{-\frac{1}{K}}^\infty h(y)e^{-\frac{\pi i}{4} (Ky+1)n}dy .
%\end{split}\end{equation}
%Being $\emph{supp}(h)\subset (h_1,h_2)$ with $0<h_1<h_2$ fixed, we can use $Ky\pm1 \gg K$ to show that the above is
%$$ \ll  \sum_{n\not=0} \int_{-\frac{1}{K}}^\infty h(y)e^{-\frac{\pi i}{4} (Ky+1)n}dy \ll \frac{1}{K^A} $$
%by integrating by parts sufficiently many times. 
%}
\end{re}

By Mellin inversion 
$$ \psi(x)=\frac{1}{2\pi i}\int_{(\sigma)}\tilde\psi(z)x^{-z}dz $$
for some $\frac{1}{2}<\sigma<1$, using \eqref{5apr.11} we have 
\begin{equation}\label{M_A(l,X)}
M_A(l,X)  
=  \frac{4}{HK}\sum_{k\equiv 0\;(4)} h\left(\frac{k-1}{K}\right)\frac{1}{2\pi i}\int_{(\sigma)}\tilde\psi(z)X^z \sum_{m=1}^\infty F_A(m,\tfrac{1}{2}+z) \sum_{f\in H_k}\!\! {\vphantom{\sum}}^{h}  \lambda(m) \lambda(l) dz.
\end{equation}
We are now ready to apply Petersson trace formula (see \cite[Proposition 14.5]{Iwaniec-Kowalski})
\begin{equation}\label{PeterssonFormula}
 \sum_{f\in H_k}\!\! {\vphantom{\sum}}^{h} \lambda(m)\lambda(n) = \delta_{m,n}+2\pi i^{-k}\sum_c\frac{S(m,n;c)}{c}J_{k-1}\left(\frac{4\pi\sqrt{mn}}{c}\right)
\end{equation}
where 
\begin{equation}\label{KloostermanSum}
S(m,n;c)=\sum_{a=1}^{c} {\vphantom{\sum}}'  e\left(\frac{am+\overline{a}n}{c}\right)
\end{equation}
is the usual Kloosterman sum and $J_{k-1}$ denotes the usual $J$-Bessel function of order $k-1$.
An application of \eqref{PeterssonFormula} to \eqref{M_A(l,X)} then yields
$$ M_A(l,X) = \mathcal D_A(l,X) + \mathcal Kl_A(l,X) $$
with
\begin{equation}\begin{split}\notag
\mathcal D_A(l,X) &=  \frac{1}{2\pi i}\int_{(\sigma)}\tilde\psi(z)X^z  F_A(l,\tfrac{1}{2}+z) \frac{4}{HK}\sum_{k\equiv 0\;(4)} h\left(\frac{k-1}{K}\right) dz\\
%&= \frac{1}{2\pi i}\int_{(\sigma)}\tilde\psi(z)X^z  F_A(l,\tfrac{1}{2}+z) \Big[1+O(K^{-A})\Big] dz
&= \frac{1}{2\pi i}\int_{(\sigma)}\tilde\psi(z)X^z  F_A(l,\tfrac{1}{2}+z)  dz + O(K^{-A})
\end{split}\end{equation}
by Remark \ref{RemTotalMeasure} 
and
\begin{equation}\begin{split}\notag
\mathcal Kl_A(l,X) & = \frac{2\pi}{2\pi i}\int_{(\sigma)}\tilde\psi(z) X^z \sum_{m=1}^{\infty}F_A(m,\tfrac{1}{2}+z)\sum_{c=1}^{\infty}\frac{S(l,m;c)}{c}\\
&\hspace{2cm}\times\frac{4}{HK}\sum_{k\equiv 0\;(4)} h\left(\frac{k-1}{K}\right) J_{k-1}\bigg(\frac{4\pi\sqrt{ml}}{c}\bigg)dz.
\end{split}\end{equation}

By \eqref{5apr.12}, one easily sees that the diagonal term matches with what the recipe predicts %.
(see Conjecture 1.5.5, Lemma 3.1.2 and Equation (4.5.4) in \cite{CFKRS}); in particular $\mathcal D_A(l,X)$ corresponds to the 0-swap term in \eqref{7june.1}.

\section{The Kloosterman part}\label{S5}

Let's consider $\mathcal Kl_A(l,X)$, from which we want to extract the 1-swap terms. As mentioned above, here we take advantage of the average over the weight $k$, that makes the behavior of the Bessel functions much easier. In particular, from \cite[Proposition 8.1]{ILS}, for $x>0$, $K\geq 1$, and $h$ a nice test function, we have
\begin{equation}\begin{split}\label{20mar.3}
4\sum_{k\equiv 0\;(4)} h\left(\frac{k-1}{K}\right)J_{k-1}(x) 
&= h\left(\frac{x}{K}\right) -\frac{K}{\sqrt{x}}\Im\left(e^{i(x-\frac{\pi}{4})}\hbar\left(\frac{K^2}{2x}\right)\right) + O\left(\frac{x}{K^3}\right) 
\end{split}\end{equation}
where 
\begin{equation}\label{20mar.4}
\hbar(v):=\int_0^{\infty}\frac{h(\sqrt{u})}{\sqrt{2\pi u}} e^{iuv} du .
\end{equation}

First, we truncate the sum over $m$ in $\mathcal Kl_A(l,X)$ at a negligible cost. Indeed, by a contour shift to the right, one can prove that, for any $\varepsilon>0$, the contribution from $m>X^{1+\varepsilon}$ to $\mathcal Kl_A(l,X)$ is $\ll K^{-1+\varepsilon}X^{-A}$ for any $A>0$.

Therefore
\begin{equation}\begin{split}\notag
\mathcal Kl_A(l,X) & = \frac{2\pi}{2\pi i}\int_{(\sigma)}\tilde\psi(z) X^z \sum_{m\ll X^{1+\varepsilon}}F_A(m,\tfrac{1}{2}+z)\sum_{c=1}^{\infty}\frac{S(l,m;c)}{c}\\
&\hspace{2cm}\times\frac{4}{HK}\sum_{k\equiv 0\;(4)} h\left(\frac{k-1}{K}\right) J_{k-1}\bigg(\frac{4\pi\sqrt{ml}}{c}\bigg)dz + O(K^{-A}).
\end{split}\end{equation}

We now apply \eqref{20mar.3}. Note that, if $Xl=K^\eta$ with $\eta<4$, then $\frac{4\pi\sqrt{ml}}{c}\ll \frac{K^{\varepsilon}\sqrt{Xl}}{c}\ll \frac{K^{2-\varepsilon}}{c}$. Therefore, being $\hbar$ a Schwartz function, the second term on the right hand side of \eqref{20mar.3} is not significant in our application (see \cite{ILS}, Remarks after Corollary 8.2).
%\cite[Proposition 8.1]{ILS} in the form
%$$ 4\sum_{k\equiv 0\;(4)}h\left(\frac{k-1}{K}\right)J_{k-1}(x) 
%= h\left(\frac{x}{K}\right) 
%- \frac{K}{\sqrt{x}}\Im\bigg( e^{i(x-\frac{\pi}{4})}\hbar\left(\frac{K^2}{2x}\right) \bigg)
%+ O\left(\frac{x}{K^3}\right) ,$$ 
%with
%$$ \hbar(v) = \int_0^{\infty} \frac{h(\sqrt{u})}{\sqrt{2\pi u}} e^{iuv} du .$$
%Then we have
Hence we have
\begin{equation}\begin{split}\label{23mar.0}
\mathcal Kl_A(l,X)  
&= \mathcal Kl_A^*(l,X)  + O\left(\frac{Xl}{K^{4-\varepsilon}}\right)
\end{split}\end{equation}
with
%\begin{equation}\begin{split}\notag
%\mathcal Kl_A^*(l,X) := \frac{2\pi}{2\pi i}\int_{(\sigma)}\tilde\psi(z) X^z \!\!\sum_{m\ll X^{1+\varepsilon}}F_A(m,\tfrac{1}{2}+z)\!\!\sum_{c\ll X^\varepsilon\frac{\sqrt{Xl}}{K^{1-\varepsilon}}\!}\frac{S(l,m;c)}{c}
%\frac{1}{HK}h\left(\frac{4\pi\sqrt{ml}}{cK}\right) dz .
%\end{split}\end{equation}
\begin{equation}\begin{split}\notag
\mathcal Kl_A^*(l,X) := \frac{2\pi}{2\pi i}\int_{(\sigma)}\tilde\psi(z) X^z \!\!\sum_{m\ll X^{1+\varepsilon}}F_A(m,\tfrac{1}{2}+z)
\sum_{c=1}^{\infty}\frac{S(l,m;c)}{c}
\frac{1}{HK}h\left(\frac{4\pi\sqrt{ml}}{cK}\right) dz .
\end{split}\end{equation}

%\textcolor{blue}{
%The first term from the RHS of \eqref{20mar.3} is just a substitution. The second term is negligible (it contributes $\ll K^{-A}$) for what we said above 
%$$
%x = \frac{4\pi\sqrt{ml}}{c} \ll \frac{\sqrt{X^{1+\varepsilon}l}}{c} \ll \frac{X^{\varepsilon/2}\sqrt{K^\eta}}{c} \ll \frac{K^{2\varepsilon}K^{2-10\varepsilon}}{c}\ll \frac{K^{2-\varepsilon}}{c} 
%\Longrightarrow\frac{K^2}{x} \gg \frac{cK^2}{K^{2-\varepsilon}} \gg cK^\varepsilon .$$
%The error term is (we use $|S(l,m,n)|\leq \tau(c)\sqrt{(l,m,c)}\sqrt{c}\ll \sqrt{l}\;c^{1/2+\varepsilon}$)
%\begin{equation}\begin{split}\notag
%\mathcal E &\ll X^\sigma\sum_{m\ll X^{1+\varepsilon}}\frac{\tau_A(m)}{m^{\frac{1}{2}+\sigma}} \sum_{c=1}^{\infty}\frac{|S(l,m;c)|}{c}\frac{\sqrt{ml}}{K^4c}
%\ll \frac{X^\sigma\sqrt{l}}{K^4}\sum_{m\ll X^{1+\varepsilon}}\frac{1}{m^{\sigma-\varepsilon}} \sum_{c=1}^{\infty}\frac{|S(l,m;c)|}{c^2}\\
%&\ll\frac{X^\sigma l}{K^4}\sum_{m\ll X^{1+\varepsilon}}\frac{1}{m^{\sigma-\varepsilon}} \sum_{c=1}^{\infty}\frac{1}{c^{3/2-\varepsilon}}
%%\ll K^{\varepsilon}  \frac{X^\sigma l}{K^4} \times \sum_{m\ll X^{1+\varepsilon}}\frac{1}{m^{\sigma}} \times 1
%\ll   \frac{X^\sigma l}{K^4} \times X^{1-\sigma+\varepsilon}\times 1
%\ll \frac{Xl}{K^{4-\varepsilon}} .
%\end{split}\end{equation}
%}

By definition of $h$, we can truncate the sum over $c$ at height $\ll \frac{\sqrt{Xl}}{K^{1-\varepsilon}}$. Moreover, with a contour shift to the right, we can re-extend the sum over $m$ at a negligible cost,
%
%\textcolor{blue}{
%\begin{equation}\begin{split}\notag
%\mathcal Kl_A^{m>X}(l,X) & 
%= \frac{2\pi}{HK}\frac{1}{2\pi i}\int_{(\sigma)}\tilde\psi(z) X^z \sum_{m\gg X^{1+\varepsilon}}F_A(m,\tfrac{1}{2}+z)\sum_{c\ll \frac{\sqrt{Xl}}{K^{1-\varepsilon}}}\frac{S(l,m;c)}{c} h\bigg(\frac{4\pi\sqrt{ml}}{cK}\bigg)dz\\
%&= \frac{2\pi}{HK}\frac{1}{2\pi i}\int_{(\frac{AK}{\varepsilon})}\tilde\psi(z) X^z \sum_{m\gg X^{1+\varepsilon}}F_A(m,\tfrac{1}{2}+z)\sum_{c\ll \frac{\sqrt{Xl}}{K^{1-\varepsilon}}}\frac{S(l,m;c)}{c} h\bigg(\frac{4\pi\sqrt{ml}}{cK}\bigg)dz\\
%&\ll \frac{X^\frac{AK}{\varepsilon}}{K}\sum_{m\gg X^{1+\varepsilon}}\frac{\tau_A(m)}{m^{\frac{1}{2}+\frac{AK}{\varepsilon}}} \sum_{c\ll \frac{\sqrt{Xl}}{K^{1-\varepsilon}}}\frac{\sqrt{l}\; c^{1/2+\varepsilon}}{c}
%\ll K^{\varepsilon} \frac{X^{\frac{AK}{\varepsilon}}\sqrt{l}}{K}\sum_{m\gg X^{1+\varepsilon}}\frac{1}{m^{\frac{1}{2}+\frac{AK}{\varepsilon}-\varepsilon}} \sum_{c\ll \frac{\sqrt{Xl}}{K^{1-\varepsilon}}}\frac{1}{\sqrt{c}}\\
%&\ll K^{\varepsilon} \frac{X^{\frac{AK}{\varepsilon}}\sqrt{l}}{K}\frac{1}{(X^{1+\varepsilon})^{-\frac{1}{2}+\frac{AK}{\varepsilon}-\varepsilon}} \sqrt{K}
%\ll {\frac{\sqrt{l}K^\varepsilon}{K}} \frac{X^{\frac{AK}{\varepsilon}}}{X^{-\frac{1}{2}+\frac{AK}{\varepsilon}-\varepsilon-\frac{\varepsilon}{2}+AK-\varepsilon^2}}\\
%&\ll \sqrt{\frac{{Xl}}{K^{1-\varepsilon}}} \frac{1}{X^{AK}}.
%\end{split}\end{equation}
%}
%
getting
\begin{equation}\begin{split}\notag
\mathcal Kl_A(l,X)  
&= \mathcal Kl_A^{**}(l,X)  + O\left(\frac{Xl}{K^{4-\varepsilon}}\right)
\end{split}\end{equation}
with
\begin{equation}\begin{split}\notag
\mathcal Kl_A^{**}(l,X) := \frac{2\pi}{HK}\frac{1}{2\pi i}\int_{(\sigma)}\tilde\psi(z) X^z \sum_{m=1}^{\infty} F_A(m,\tfrac{1}{2}+z)\sum_{c\ll \frac{\sqrt{Xl}}{K^{1-\varepsilon}}}\frac{S(l,m;c)}{c} h\left(\frac{4\pi \sqrt{ml}}{cK}\right) dz.
\end{split}\end{equation}

%We immediately see that 
%\begin{equation}\label{23mar.1} \mathcal Kl_A^2(l,X) \ll K^{-A} \end{equation}
%for any $A>0$, since for $c\geq 1$, $m\leq X^{1+\varepsilon} \;\forall \varepsilon>0$, $Xl\ll K^{4-\delta}$ for some fixed $\delta>0$, we have
%$$\frac{4\pi\sqrt{ml}}{c}\leq\frac{15\sqrt{X^{1+\varepsilon}l}}{1}\leq 15X^{\frac{\varepsilon}{2}}\sqrt{K^{4-\delta}}\leq 15K^{2-\frac{\delta}{2}+2\varepsilon}$$
%i.e.
%$$ \frac{K^2c}{8\pi\sqrt{ml}} \geq \frac{1}{30}K^{\frac{\delta}{2}-2\varepsilon} \gg K^{\delta/3} $$
%taking $\varepsilon$ sufficiently small with respect to $\delta$. Then, integrating by parts sufficiently many times, we have
%$$ \hbar\left(\frac{K^2c}{8\pi\sqrt{ml}}\right) \ll K^{-B} $$
%for any $B>0$. Therefore \eqref{23mar.1} follows and, together with \eqref{23mar.0} implies
%\begin{equation}\begin{split}\label{23mar.2}
%\mathcal Kl_A(l,X)  
%&= \mathcal Kl_A^1(l,X) + O\left(\frac{Xl}{K^{4-\varepsilon}}\right)
%\end{split}\end{equation}

By Mellin inversion 
$$ h(x) = \frac{1}{2\pi i}\int_{(\nu)}\tilde h(w)x^{-w}dw = \frac{2}{2\pi i}\int_{(\frac{\nu}{2})}\tilde h(2w)x^{-2w}dw $$
for some $\nu>2$, opening the Kloosterman sum, we write the above as
\begin{equation}\begin{split}\notag
\mathcal Kl_A^{**}(l,X)% :&= \frac{2\pi}{HK}\frac{1}{2\pi i}\int_{(\sigma)}\tilde\psi(z) X^z 
%\frac{1}{2\pi i}\int_{(\nu)} \tilde h(w)\left(\frac{4\pi\sqrt{l}}{K}\right)^{-w} \\
%&\hspace{2cm}\times \sum_{m=1}^{\infty}\frac{F_A(m,\tfrac{1}{2}+z)}{m^{w/2}}%\sum_{c\ll X^\varepsilon\frac{\sqrt{Xl}}{K^{1-\varepsilon}}}\frac{S(l,m;c)}{c^{1-w}}
%\;dw\; dz \\
%&= \frac{1}{HK}\frac{1}{2\pi i}\int_{(\sigma)}\tilde\psi(z) X^z 
%\frac{1}{2\pi i}\int_{(\nu)} \tilde h(2w) \frac{K^{2w}(4\pi)^{1-2w}}{l^{w}}  \\
%&\hspace{2cm}\times \sum_{m=1}^{\infty} \frac{F_A(m,\tfrac{1}{2}+z)}{m^{w}}\sum_{c\ll X^\varepsilon\frac{\sqrt{Xl}}{K^{1-\varepsilon}}}\frac{S(l,m;c)}{c^{1-2w}}
%\;dw\; dz\\
&= \frac{1}{HK}\frac{1}{2\pi i}\int_{(\sigma)}\tilde\psi(z) X^z 
\frac{1}{2\pi i}\int_{(\frac{\nu}{2})} \tilde h(2w) \frac{K^{2w}(4\pi)^{1-2w}}{l^{w}}  \\
&\hspace{2cm}\times
\sum_{c\ll \frac{\sqrt{Xl}}{K^{1-\varepsilon}}}\frac{1}{c^{1-2w}} \sum_{a=1}^{c} {\vphantom{\sum}}'  e\left(\frac{\overline{a}l}{c}\right) 
\sum_{m=1}^{\infty} \frac{F_A(m,\tfrac{1}{2}+z)e(\frac{am}{c})}{m^{w}}
\;dw\; dz .
\end{split}\end{equation}
%
%\textcolor{red}{\\ **** \\}
%
%Now we handle the remaining sums (over $m$ and $c$) in the same way as in the previous version, as the second average does not change the arithmetic of the problem.

\begin{re}\label{RemarkEstermann}
Note that the inner sum is a (generalized) \lq\lq Estermann function.\rq\rq \; The classical Estermann function has been studied in literature \cite{Estermann, Conrey2/5}; typically one writes the additive character in terms of Dirichlet characters and expect poles for $\chi=\chi_0$ and entire $L$-functions otherwise. The same phenomenon appears also in this generalized case, as shown by the following considerations together with Lemma \ref{lemmaEF} and Lemma \ref{lemmaEFrefined}.
\end{re}

To carry on our analysis of the Kloosterman part from Petersson formula with the aim of detecting the 1-swap terms from the recipe, we express the additive character $e(\frac{am}{c})$ appearing in $\mathcal Kl_A(l,X)$ in terms of Dirichlet characters and Gauss sums. Namely (recall that $(a,c)=1$):
\begin{equation}\label{6apr.3} 
e\left(\frac{am}{c}\right) 
= e\left(\frac{am/(c,m)}{c/(c,m)}\right) 
=\frac{1}{\varphi(c/(m,c))}\sum_{\chi\;(c/(m,c))}\tau(\overline\chi)\chi\left(\frac{am}{(c,m)}\right), 
\end{equation}
where
$ \tau(\chi) = \sum_{n\;(\text{mod }q)}\chi(n)e(\tfrac{n}{q}) $
is the usual Gauss sum associated with a character $\chi$ modulo $q$. 
By \eqref{6apr.3} we get (calling $(m,c)=g$ and then performing the change of variable $m\to mg$)
\begin{equation}\begin{split}\notag
%\mathcal Kl_A^{**}(l,X) &=  
%\frac{1}{HK}\frac{1}{2\pi i}\int_{(\sigma)}\tilde\psi(z) X^z 
%\frac{1}{2\pi i}\int_{(\frac{\nu}{2})} \tilde h(2w) \frac{K^{2w}(4\pi)^{1-2w}}{l^{w}} \\
%&\quad \times \sum_{c\ll X^\varepsilon \frac{\sqrt{Xl}}{K^{1-\varepsilon}}} \frac{1}{c^{1-2w}}\sum_{a=1}^{c} {\vphantom{\sum}}'  e\left(\frac{\overline{a}l}{c}\right) 
%\sum_{g|c}\frac{g^{-w}}{\varphi(\frac{c}{g})}\sum_{\chi\;(\frac{c}{g})}\tau(\overline\chi) \chi(a) \sum_{\substack{m=1 \\ (m,\frac{c}{g})=1}}^\infty \frac{F_A(mg,\tfrac{1}{2}+z)\chi(m) }{m^w}  dw\;dz\\
&=\frac{1}{HK}\frac{1}{2\pi i}\int_{(\sigma)}\tilde\psi(z) X^z 
\frac{1}{2\pi i}\int_{(\frac{\nu}{2})} \tilde h(2w) \frac{K^{2w}(4\pi)^{1-2w}}{l^{w}} \\
&\quad \times \sum_{c\ll \frac{\sqrt{Xl}}{K^{1-\varepsilon}}} \frac{1}{c^{1-2w}}\sum_{a=1}^{c} {\vphantom{\sum}}'  e\left(\frac{\overline{a}l}{c}\right) 
\sum_{g|c}\frac{g^{-w}}{\varphi(\frac{c}{g})}\sum_{\chi\;(\frac{c}{g})}\tau(\overline\chi) \chi(a) \mathcal F_A(w,\tfrac{1}{2}+z,\chi,c,g)  dw\;dz
\end{split}\end{equation}
since $\Re(\frac{1}{2}+z)=\frac{1}{2}+\sigma>1$ and $\Re(w)>1$. In view of the meromorphic continuation of $\mathcal F_A(w,\tfrac{1}{2}+z,\chi,c,g)$ in the region $\Re(\frac{1}{2}+z+w)>\frac{1}{2}$, given by Lemma \ref{lemmaEFrefined}, we can move the $w$-integral to the left, on the line $(\nu')$, with $0<\sigma+\nu'<\frac{1}{2}$. An admissible choice is $\sigma=\frac{1}{2}+\varepsilon$, $\nu=3$, $\nu'=-\frac{1}{2}$.
Doing so, we encounter poles at $w+z = \frac{1}{2}-\alpha$, $\alpha\in A$ when $\chi=\chi_0$ (see Lemma \ref{lemmaEF} and Lemma \ref{lemmaEFrefined}) and we get
\begin{equation}\label{21mar.1}
\mathcal Kl_A^{**}(l,X) = \sum_{\alpha\in A}\mathcal P_\alpha +\mathcal I,
\end{equation}
where
\begin{equation}\begin{split}\label{21mar.2}
\mathcal I 
&=\frac{1}{HK}\frac{1}{2\pi i}\int_{(\sigma)}\tilde\psi(z) X^z 
\frac{1}{2\pi i}\int_{(\nu')} \tilde h(2w) \frac{K^{2w}(4\pi)^{1-2w}}{l^{w}} \\
&\quad \times \sum_{c\ll \frac{\sqrt{Xl}}{K^{1-\varepsilon}}} \frac{1}{c^{1-2w}}\sum_{a=1}^{c} {\vphantom{\sum}}'  e\left(\frac{\overline{a}l}{c}\right) 
\sum_{g|c}\frac{g^{-w}}{\varphi(\frac{c}{g})}\sum_{\chi\;(\frac{c}{g})}\tau(\overline\chi) \chi(a) \mathcal F_A(w,\tfrac{1}{2}+z,\chi,c,g)  dw\;dz
\end{split}\end{equation}
and
\begin{equation}\begin{split}\label{21mar.3}
\mathcal P_\alpha &= \frac{1}{HK}\frac{1}{2\pi i}\int_{(\sigma)}\tilde\psi(z) X^z 
\tilde h(1-2z-2\alpha) K^{1-2z-2\alpha}(4\pi)^{2z+2\alpha}\\
&\quad \times \res_{w=\frac{1}{2}-z-\alpha}\bigg( \sum_{c\ll \frac{\sqrt{Xl}}{K^{1-\varepsilon}}} \frac{1}{c^{1-2w}}\frac{1}{l^w}\sum_{a=1}^{c} {\vphantom{\sum}}'  e\left(\frac{\overline{a}l}{c}\right) 
\sum_{g|c}\frac{g^{-w}}{\varphi(\frac{c}{g})}\tau(\overline\chi_0) \mathcal F_A(w,\tfrac{1}{2}+z,\chi_0,c,g)  \bigg)dz.
\end{split}\end{equation}

First we bound $\mathcal I$, assuming the Generalized Lindel\"of Hypothesis for Dirichlet $L$-functions. Let's fix $\sigma=\frac{1}{2}+\varepsilon$ and $\nu'=-\frac{1}{2}$. Under GLH, $\mathcal F_A(w,\frac{1}{2}+z,\chi,c,g)\ll (Kc)^{r\varepsilon}$ for $\Re(w+\frac{1}{2}+z)=\frac{1}{2}+\varepsilon$. Therefore we can bound
\begin{equation}\notag
\mathcal I 
\ll \frac{1}{K} X^{\frac{1}{2}+\varepsilon}  \frac{K^{-1}}{l^{-1/2}}\sum_{c\ll \frac{\sqrt{Xl}}{K^{1-\varepsilon}}} \frac{1}{c^2} 
\sum_{g|c}\frac{g^{1/2}}{\varphi(\frac{c}{g})}\sum_{\chi\;(\frac{c}{g})}|\tau(\overline\chi)| \bigg|\sum_{a=1}^{c} {\vphantom{\sum}}'  e\left(\frac{\overline{a}l}{c}\right) 
 \chi(a)\bigg| (Kc)^{r\varepsilon}.
\end{equation}
Moreover, for $\chi$ a character modulo $\tfrac{c}{g}$, we recall that $|\tau(\overline\chi)| \leq \sqrt{c/g}$. Also, a character $\chi$ (mod $\frac{c}{g}$) induces a character $\chi_1=\chi\chi_0$ (mod $c$), where $\chi_0$ is the principal character (mod $c$); then we have 
\begin{equation}\begin{split}\notag
\bigg|\sum_{a=1}^{c} {\vphantom{\sum}}'  e\left(\frac{al}{c}\right)\chi(a) \bigg|
&%= \sum_{a=1}^{c}  e\left(\frac{al}{c}\right)\chi(a)\chi_0(a) 
= \bigg|\sum_{a=1}^{c}   e\left(\frac{al}{c}\right)\chi_1(a) \bigg|%\\&\
\leq \sqrt{c}\sum_{d|(l,c)}d
\end{split}\end{equation}
since, if $\chi_1$ modulo $c$ is induced by the primitive character $\chi^*$ modulo $h$, with $h|c$, then by \cite[Lemma 7.1]{Petrow-Young}
\begin{equation}\begin{split}\notag
\sum_{a=1}^{c}   e\left(\frac{al}{c}\right)\chi_1(a)
= \tau(\chi^*) \sum_{d|(l,\frac{c}{h})} d \;\overline{\chi^*}\bigg(\frac{l}{d}\bigg)\chi^*\bigg(\frac{c}{dh}\bigg)\mu\bigg(\frac{c}{dh}\bigg).
\end{split}\end{equation}
Combining these estimates and bounding trivially the remaining sums, we get (under GLH)
\begin{equation}\begin{split}\notag
\mathcal I 
%&\ll K^{r\varepsilon} \frac{\sqrt{Xl}}{K^2} \sum_{c\ll  \frac{\sqrt{Xl}}{K^{1-\varepsilon}}} \frac{c^{r\varepsilon}}{c^2} 
%\sum_{g|c}\frac{\sqrt{g}}{\varphi(\frac{c}{g})}\sum_{\chi\;(\frac{c}{g})}\sqrt{\frac{c}{g}} \bigg|\tau(\chi_1)\sum_{d|(l,g)} d\;\overline\chi_1^*\bigg(\frac{l}{d}\bigg)\mu\bigg(\frac{g}{d}\bigg)\bigg|\\
%&\ll K^{r\varepsilon} \frac{\sqrt{Xl}}{K^2} \sum_{c\ll \frac{\sqrt{Xl}}{K^{1-\varepsilon}}} \frac{c^{r\varepsilon}}{c^2} 
%\sum_{g|c}\frac{\sqrt{g}}{\varphi(\frac{c}{g})}\sum_{\chi\;(\frac{c}{g})}\sqrt{\frac{c}{g}}\sqrt{c} \sum_{d|(l,g)} d \\
&\ll K^{r\varepsilon} \frac{\sqrt{Xl}}{K^2} \sum_{c\ll \frac{\sqrt{Xl}}{K^{1-\varepsilon}}} \frac{c^{r\varepsilon}}{c} 
\sum_{g|c}\frac{1}{\varphi(\frac{c}{g})}\sum_{\chi\;(\frac{c}{g})} \sum_{d|(l,c)} d 
\ll K^{r\varepsilon} \frac{\sqrt{Xl}}{K^2} \sum_{c\ll \frac{\sqrt{Xl}}{K^{1-\varepsilon}}} \frac{c^{r\varepsilon}}{c}  \sum_{d|(l,c)} d \\
%&\ll K^{r\varepsilon} \frac{\sqrt{Xl}}{K^2}  \sum_{d|l}d \sum_{\substack{c\ll X^\varepsilon \frac{\sqrt{Xl}}{K^{1-\varepsilon}}\\ c\equiv 0\;(d)}} \frac{c^{r\varepsilon}}{c}  
%&\ll K^{r\varepsilon} \frac{\sqrt{Xl}}{K^2}  \sum_{d|l}d \sum_{c\ll X^\varepsilon \frac{\sqrt{Xl}}{dK^{1-\varepsilon}}} \frac{c^{r\varepsilon}}{cd} 
&\ll K^{r\varepsilon} \frac{\sqrt{Xl}}{K^2}  \sum_{d|l} \sum_{c\ll  \frac{\sqrt{Xl}}{dK^{1-\varepsilon}}} \frac{c^{r\varepsilon}}{c} 
\ll K^{2r\varepsilon} \frac{\sqrt{Xl}}{K^2}.
\end{split}\end{equation}

\section{The contribution from the poles}\label{S6}
Finally, for any $\tilde\alpha\in A$, we show that $\mathcal P_{\tilde\alpha}$ defined in \eqref{21mar.3} corresponds to the 1-swap $\tilde\alpha\rightarrow -\tilde\alpha$ term from the recipe, i.e. the one corresponding to the term $V=\{\tilde\alpha\}$ in the notations of Theorem \ref{thm}.
We recall that $\sum_{a=1}^{c} {\vphantom{\sum}}'  e\left(\frac{\overline{a}l}{c}\right)  = R_c(l)$ the Ramanujan sum
and
$\tau(\overline\chi_0) = \mu(c/g)$ if $\chi_0$ is mod$\frac{c}{g}$; then
\begin{equation}\begin{split}\notag
\mathcal P_{\tilde\alpha} %&= \frac{1}{HK}\frac{1}{2\pi i}\int_{(\sigma)}\tilde\psi(z) X^z 
%\tilde h(1-2z-2\alpha) \frac{K^{1-2z-2\alpha}(4\pi)^{2z+2\alpha}}{l^{\frac{1}{2}-z-\alpha}}\\
%&\quad \times \res_{w=\frac{1}{2}-z-\tilde\alpha}\bigg( \sum_{c\ll X^\varepsilon \frac{\sqrt{Xl}}{K^{1-\varepsilon}}} \frac{1}{c^{1-2w}}\sum_{a=1}^{c} {\vphantom{\sum}}'  e\left(\frac{\overline{a}l}{c}\right) 
%\sum_{g|c}\frac{g^{-w}}{\varphi(\frac{c}{g})}\tau(\overline\chi_0) \mathcal F_A(w,\tfrac{1}{2}+z,\chi_0,c,g)  \bigg)dz\\
&=\frac{1}{HK}\frac{1}{2\pi i}\int_{(\sigma)}\tilde\psi(z) X^z 
\tilde h(1-2z-2\alpha) K^{1-2z-2\alpha}(4\pi)^{2z+2\alpha}\\
&\quad \times \res_{w=\frac{1}{2}-z-\tilde\alpha}\bigg( \sum_{c\ll \frac{\sqrt{Xl}}{K^{1-\varepsilon}}} \frac{R_c(l)}{c^{1-2w}}\frac{1}{l^w}
\sum_{g|c}\frac{g^{-w}\mu(\frac{c}{g})}{\varphi(\frac{c}{g})}\mathcal F_A(w,\tfrac{1}{2}+z,\chi_0,c,g)  \bigg)dz\\
&=\frac{1}{2\pi i}\int_{(\sigma)}\tilde\psi(z) X^z 
\frac{\tilde h(1-2z-2\alpha)}{H} \left(\frac{K^2}{(4\pi)^2}\right)^{-z-\alpha}\\
&\hspace{2cm} \times \res_{w=\frac{1}{2}-z-\tilde\alpha}\bigg( \sum_{c\ll \frac{\sqrt{Xl}}{K^{1-\varepsilon}}} \frac{R_c(l)}{c^{1-2w}}\frac{1}{l^w}
\sum_{m=1}^{\infty} \frac{\mu(\frac{c}{(m,c)})}{\varphi(\frac{c}{(m,c)})} \frac{F_A(m,\frac{1}{2}+z)}{m^w}  \bigg)dz.
\end{split}\end{equation}

%\textcolor{blue}{
%Being
%\begin{equation}\begin{split}\notag
%\sum_{m=1}^{\infty} \frac{\mu(\frac{c}{(m,c)})}{\varphi(\frac{c}{(m,c)})}&\frac{F_A(m,\frac{1}{2}+z)}{m^w}
%= \sum_{\substack{m=1 \\ m\equiv 0\;(g)}}^{\infty} \frac{\mu(\frac{c}{g})}{\varphi(\frac{c}{g})}\frac{F_A(m,\frac{1}{2}+z)}{m^w}\\
%&= \sum_{g|c} \frac{\mu(\frac{c}{g})}{\varphi(\frac{c}{g})} \sum_{\substack{m=1 \\ (m,c)=g}}^{\infty} \frac{F_A(m,\frac{1}{2}+z)}{m^w}
%= \sum_{g|c} \frac{g^{-w}\mu(\frac{c}{g})}{\varphi(\frac{c}{g})} \sum_{\substack{m=1 \\ (m,c)=g}}^{\infty} \frac{F_A(mg,\frac{1}{2}+z)}{m^w}.
%\end{split}\end{equation}
%}

We note that, at the cost of a negligible error term, we can extend the sum over $c$ to infinity. 
To justify this, one uses the Perron formula to express the truncated sum over $c$ as an integral; after a contour shift, the pole at 0 gives the completed sum, while the contributions from the remaining integral is as small as $K^{-\varepsilon}$. 
Then we look at
\begin{equation}\begin{split}\label{24mar.11}\notag
%\mathcal P_{\tilde\alpha}= 
&\frac{1}{2\pi i}\int_{(\sigma)}\tilde\psi(z) X^z 
\frac{\tilde h(1-2z-2\alpha)}{H} \left(\frac{K^2}{(4\pi)^2}\right)^{-z-\alpha}\\
&\hspace{2cm} \times \res_{w=\frac{1}{2}-z-\tilde\alpha}\bigg( \sum_{c=1}^{\infty} \frac{R_c(l)}{c^{1-2w}}\frac{1}{l^w}
\sum_{m=1}^{\infty} \frac{\mu(\frac{c}{(m,c)})}{\varphi(\frac{c}{(m,c)})} \frac{F_A(m,\frac{1}{2}+z)}{m^w}  \bigg)dz
\end{split}\end{equation}
and show that it corresponds to the 1-swap term $\tilde\alpha \to -\tilde\alpha$. By the Poisson summation formula, for any bounded $\gamma$, one can easily see that
$$ \frac{\tilde h(1-2\gamma)}{H}\left(\frac{K^2}{(4\pi)^2}\right)^{-\gamma}
= \frac{4}{HK}\sum_{k\equiv0\;(4)}h\left(\frac{k-1}{K}\right)\left(\frac{k}{4\pi}\right)^{-2\gamma} + O\bigg(\frac{1}{K}\bigg).$$
%
%\textcolor{blue}{
%indeed
%\begin{equation}\begin{split}\notag
%LHS &=  \frac{\tilde h(1-2\gamma)}{H}\left(\frac{K^2}{(4\pi)^2}\right)^{-\gamma}
%= \frac{1}{H}\int_0^{\infty} h(t)t^{-2\gamma} dt \left(\frac{K^2}{(4\pi)^2}\right)^{-\gamma}= \frac{1}{H}\int_0^{\infty} h(t) \left(\frac{tK}{4\pi}\right)^{-2\gamma}dt
%\end{split}\end{equation}
%\begin{equation}\begin{split}\notag
%RHS &= \frac{4}{HK}\sum_{k\equiv0\;(4)}h\left(\frac{k-1}{K}\right)\left(\frac{k^2}{16\pi^2}\right)^{-\gamma}
%= \frac{4}{HK}\sum_{m}h\left(\frac{4m-1}{K}\right)\left(\frac{m^2}{\pi^2}\right)^{-\gamma}\\
%&= \frac{4}{HK}\sum_{n}\mathcal F\bigg(h\left(\frac{4*-1}{K}\right)\left(\frac{*^2}{\pi^2}\right)^{-\gamma}\bigg)(n)
%= \frac{4}{HK}\sum_{n}\int_0^{\infty}h\left(\frac{4x-1}{K}\right)\left(\frac{x^2}{\pi^2}\right)^{-\gamma}e^{-2\pi i x n}dx\\
%&= \frac{1}{H}\sum_{n}\int_0^{\infty}h(t)\left(\frac{Kt+1}{4\pi}\right)^{-2\gamma}e^{-2\pi i \frac{xK+1}{4} n}dt 
%= \frac{1}{H}\int_0^{\infty}h(t)\left(\frac{Kt+1}{4\pi}\right)^{-2\gamma}dt + O(K^{-A}) \\
%&= \frac{1}{H}\int_0^{\infty}h(t)\bigg(\frac{tK}{4\pi}\bigg)^{-2\gamma}\left(1+\frac{1}{tK}\right)^{-2\gamma}dt + O(K^{-A}) 
%= \frac{1}{H}\int_0^{\infty}h(t)\bigg(\frac{tK}{4\pi}\bigg)^{-2\gamma}dt + O(K^{-1}) 
%\end{split}\end{equation}
%}
%
This gives
\begin{equation}\begin{split}\label{24mar.12}\notag
\mathcal P_{\tilde\alpha}= 
&\frac{1}{2\pi i}\int_{(\sigma)}\tilde\psi(z) X^z 
\frac{4}{HK}\sum_{k\equiv0\;(4)}h\left(\frac{k-1}{K}\right) \left(\frac{k}{4\pi}\right)^{-2z-2\alpha}\\
&\hspace{2cm} \times \res_{w=\frac{1}{2}-z-\tilde\alpha}\bigg( \sum_{c=1}^{\infty} \frac{R_c(l)}{c^{1-2w}}\frac{1}{l^w}
\sum_{m=1}^{\infty} \frac{\mu(\frac{c}{(m,c)})}{\varphi(\frac{c}{(m,c)})} \frac{F_A(m,\frac{1}{2}+z)}{m^w}  \bigg)dz +O(K^{-\varepsilon}).
\end{split}\end{equation}
Therefore, we only have to prove that the residue equals $G_{A\setminus\{\tilde\alpha\}\cup\{-\tilde\alpha\}}(\tfrac{1}{2};l)$. We study the second line above by Euler product and, denoting $s=\frac{1}{2}+z$, we have
\begin{equation}\begin{split}\notag
\mathcal S(l) 
:&= \sum_{c=1}^{\infty} \frac{R_c(l)}{c^{1-2w}}\frac{1}{l^w}
\sum_{m=1}^\infty \frac{\mu(c/(m,c))}{\varphi(c/(m,c))} \frac{F_A(m,s) }{m^w} \\
&=\prod_{p^\nu || l} \frac{1}{p^{\nu w}}\sum_{c=0}^\infty \frac{R_{p^c}(l)}{p^{c(1-2w)}}\sum_{m=0}^\infty \frac{\mu(p^{c-\min(m,c)})}{\varphi(p^{c-\min(m,c)})}\frac{F_{A,p}(m,s)}{p^{mw}} .
\end{split}\end{equation}
When $p^\nu || l$, one can easily prove that 
$$ R_{p^c}(l) = \begin{cases}
1 & \text{if }c=0 \\ \varphi(p^c) & \text{if }1\leq c\leq\nu \\ -p^\nu & \text{if }c=\nu+1 \\ 0 & \text{if }c\geq\nu+2.
\end{cases} $$
Therefore, the above is
\begin{equation}\begin{split}\notag 
\mathcal S(l)=
%\prod_{p^\nu || l} \bigg(\frac{1}{p^{\nu w}}\sum_{m=0}^\infty \frac{F_{A,p}(m,s)}{p^{mw}} + \sum_{1\leq c\leq \nu} \frac{\varphi(p^c)}{p^{c(1-2w)}p^{\nu w}}\sum_{m=0}^\infty \frac{\mu(p^{c-\min(m,c)})}{\varphi(p^{c-\min(m,c)})}\frac{F_{A,p}(m,s)}{p^{mw}}  \\
%&\hspace{3cm} +\frac{-p^\nu}{p^{(\nu+1)(1-2w)}p^{\nu w}}\sum_{m=0}^\infty \frac{\mu(p^{\nu+1-\min(m,\nu+1)})}{\varphi(p^{\nu+1-\min(m,\nu+1)})}\frac{F_{A,p}(m,s)}{p^{mw}}  \bigg)\\
%=&
\prod_{p^\nu || l} \bigg(\frac{1}{p^{\nu w}}&\sum_{m=0}^\infty \frac{F_{A,p}(m,s)}{p^{mw}} -\sum_{1\leq c\leq \nu} \frac{p^{c-1}}{p^{c(1-2w)}p^{\nu w}} \frac{F_{A,p}(c-1,s)}{p^{(c-1)w}}  \\
&+ \sum_{1\leq c\leq \nu} \frac{p^c(1-\frac{1}{p})}{p^{c(1-2w)}p^{\nu w}}\sum_{m=c}^\infty \frac{F_{A,p}(m,s)}{p^{mw}} + \frac{p^\nu}{p^{(\nu+1)(1-2w)}p^{\nu w}}\frac{1}{p-1}\frac{F_{A,p}(\nu,s)}{p^{\nu w}} \\
& - \frac{p^\nu}{p^{(\nu+1)(1-2w)}p^{\nu w}}\sum_{m=\nu+1}^\infty\frac{F_{A,p}(m,s)}{p^{mw}} \bigg)
\end{split}\end{equation}
since $\mu(p^{c-\min(m,c)})=0$ if $m<c-1$ and $\varphi(p^c)=p^c(1-\frac{1}{p})$ if $c>0$. Completing the geometric series and applying \eqref{5apr.12} and \eqref{26june.1}, we get
\begin{equation}\begin{split}\notag 
\mathcal S(l)=&
\prod_{p^\nu || l}  \frac{2}{\pi}\int_0^\pi \prod_{\alpha\in A}C\left(\frac{1}{p^{s+\alpha}};\theta\right) \frac{1}{p^{\nu w}}\bigg\{\bigg[ 1+ \sum_{1\leq c\leq \nu} \frac{p^c(1-\frac{1}{p})}{p^{c(1-2w)}} - \frac{p^\nu}{p^{(\nu+1)(1-2w)}} \bigg]C\left(\frac{1}{p^w};\theta\right)  \\
&\hspace{1cm}  -\frac{1}{p^{1-2w}}\sum_{1\leq c\leq \nu}\frac{U_{c-1}(\cos \theta)}{p^{-w(c-1)}}  -  \sum_{1\leq c\leq \nu} \frac{p^c(1-\frac{1}{p})}{p^{c(1-2w)}}\sum_{m=0}^{c-1} \frac{U_m(\cos\theta)}{p^{mw}} \\ 
&\hspace{1cm} +\frac{p^\nu}{p^{(\nu+1)(1-2w)}}\frac{1}{p-1}\frac{U_\nu(\cos\theta)}{p^{\nu w}}  + \frac{p^\nu}{p^{(\nu+1)(1-2w)}}\sum_{m=0}^\nu\frac{U_m(\cos\theta)}{p^{mw}}  \bigg\}\sin^2\theta\; d\theta.
\end{split}\end{equation}

To prove that the residue of the above quantity gives the 1-swap term corresponding to $V=\{\tilde\alpha\}$, it suffices to show that% for any $\alpha\in A$ we have 
\begin{equation}\begin{split}\label{th}  
C\left(\frac{1}{p^{s+\tilde\alpha}};\theta\right) &\frac{1}{p^{\nu w}}\bigg\{\bigg[ 1+ \sum_{1\leq c\leq \nu} \frac{p^c(1-\frac{1}{p})}{p^{c(1-2w)}} - \frac{p^\nu}{p^{(\nu+1)(1-2w)}} \bigg]C\left(\frac{1}{p^w};\theta\right)  \\
&-\frac{1}{p^{1-2w}}\sum_{1\leq c\leq \nu}\frac{U_{c-1}(\cos \theta)}{p^{-w(c-1)}}  -  \sum_{1\leq c\leq \nu} \frac{p^c(1-\frac{1}{p})}{p^{c(1-2w)}}\sum_{m=0}^{c-1} \frac{U_m(\cos\theta)}{p^{mw}} \\ 
&+\frac{p^\nu}{p^{(\nu+1)(1-2w)}}\frac{1}{p-1}\frac{U_\nu(\cos\theta)}{p^{\nu w}}  + \frac{p^\nu}{p^{(\nu+1)(1-2w)}}\sum_{m=0}^\nu\frac{U_m(\cos\theta)}{p^{mw}}  \bigg\}\bigg|_{\substack{s=1/2+z \\ w=1/2-\tilde\alpha}}\\ 
= &\;C\left(\frac{1}{p^{s-\tilde\alpha}};\theta\right)\left(1-\frac{1}{p}\right)^{-1}U_\nu(\cos\theta)\bigg|_{s=1/2-z}.
\end{split}\end{equation}

\subsection{Proof of \eqref{th}}
First of all we notice that, being the set of shifts arbitrary, it suffices to show the claim for $z=0$. We also introduce some notations to make the above more readable, denoting $x=1/p^{1/2}$ and $y=1/p^{\tilde\alpha}$. Moreover, we will write only $C(X)$ for $C(X;\theta)$ and $U_m$ for the Chebyshev polynomial $U_m(\cos\theta)$, keeping in mind that these are not \lq\lq independent\rq\rq variables, as the Chebyshev polynomials satisfy (and in fact they can be defined by) the recurrence relation
\begin{equation}\label{recurrence} U_1U_M = U_{M-1} + U_{M+1} .\end{equation}
Therefore our \eqref{th} amounts to proving
\begin{equation}\begin{split}\label{th1}
&\left(\frac{x}{y}\right)^{\nu}\bigg\{\bigg[ 1+ \left(1-x^2\right)\sum_{1\leq c\leq \nu} \left(\frac{y}{x}\right)^{2c} - \frac{y^{2(\nu+1)}}{x^{2\nu}} \bigg] \\
&+C\left(\frac{x}{y}\right)^{-1}\bigg[-y^2\sum_{1\leq c\leq \nu}\left(\frac{y}{x}\right)^{c-1}U_{c-1} - \left(1-x^2\right) \sum_{1\leq c\leq \nu} \left(\frac{y}{x}\right)^{2c}\sum_{m=0}^{c-1} \left(\frac{x}{y}\right)^m U_m \\ 
&+\frac{y^{2(\nu+1)}}{x^{2\nu}}\frac{x^2}{1-x^2}\left(\frac{x}{y}\right)^{\nu}U_\nu  + \frac{y^{2(\nu+1)}}{x^{2\nu}}\sum_{m=0}^\nu\left(\frac{x}{y}\right)^m U_m\bigg]  \bigg\}\\ 
= &\;C\left(xy\right)^{-1}\left(1-x^2\right)^{-1}U_\nu.
\end{split}\end{equation}
Hence, recalling that $C(X)^{-1}=1-U_1X+X^2$, the right hand side of \eqref{th1} can be written as
\begin{equation}\notag
RHS(\nu) 
= \frac{1+x^2y^2}{1-x^2}U_\nu - \frac{xy}{1-x^2}U_1U_\nu . 
\end{equation}
Therefore, by \eqref{recurrence}, we have
\begin{equation}\begin{split}\notag
RHS(0) &= \frac{1+x^2y^2}{1-x^2} - \frac{xy}{1-x^2}U_1 , 
\quad RHS(1) = \frac{1+x^2y^2}{1-x^2}U_1 - \frac{xy}{1-x^2} U_1^2, \\
RHS(\nu) &= \frac{1+x^2y^2}{1-x^2}U_\nu - \frac{xy}{1-x^2}U_{\nu-1}  - \frac{xy}{1-x^2}U_{\nu+1} \quad \text{ for } \nu\geq 2.
\end{split}\end{equation}
The left hand side of \eqref{th1} (say $LHS(\nu)$) instead reads
%\begin{equation}\begin{split}\notag
%LHS(\nu) :=&\left(\frac{x}{y}\right)^\nu \bigg\{\bigg[ 1+ (1-x^2)\sum_{1\leq c\leq \nu} \left(\frac{y}{x}\right)^{2c} - \frac{y^{2(\nu+1)}}{x^{2\nu}} \bigg] \\
%&+C\left(\frac{x}{y}\right)^{-1}\bigg[-y^2 \sum_{1\leq c\leq \nu}\left(\frac{y}{x}\right)^{c-1}U_{c-1} - (1-x^2) \sum_{1\leq c\leq \nu} \left(\frac{y}{x}\right)^{2c}\sum_{m=0}^{c-1} \left(\frac{x}{y}\right)^{m}U_m \\ 
%&+\frac{y^{2(\nu+1)}}{x^{2\nu}}\frac{x^2}{1-x^2}\left(\frac{x}{y}\right)^{\nu}U_\nu  + \frac{y^{2(\nu+1)}}{x^{2\nu}}\sum_{m=0}^\nu\left(\frac{x}{y}\right)^{m}U_m\bigg]  \bigg\}
%\end{split}\end{equation}
%which we will write as
\begin{equation}\notag
LHS(\nu) =\left(\frac{x}{y}\right)^\nu \bigg\{ A_{\nu} +\left(1+\frac{x^2}{y^2}\right) \sum_{i=0}^\nu A_{i,\nu} U_i -\frac{x}{y}\sum_{i=0}^\nu A_{i,\nu}U_1 U_i\bigg\}
\end{equation}
with the notations
\begin{equation}\begin{split}\notag
& A_{\nu} := 1+ (1-x^2)\sum_{1\leq c\leq \nu} \left(\frac{y}{x}\right)^{2c} - \frac{y^{2(\nu+1)}}{x^{2\nu}}\\
& A_{i,\nu} :=  -y^2 \left(\frac{y}{x}\right)^i - (1-x^2) \left(\frac{x}{y}\right)^i \sum_{i+1\leq c\leq \nu} \left(\frac{y}{x}\right)^{2c} + \frac{y^{2(\nu+1)}}{x^{2\nu}} \left(\frac{x}{y}\right)^{i} \quad \text{ for } 0\leq i \leq \nu-1\\
& A_{\nu,\nu} := \frac{y^{2(\nu+1)}}{x^{2\nu}}\frac{x^2}{1-x^2}\left(\frac{x}{y}\right)^{\nu} + \frac{y^{2(\nu+1)}}{x^{2\nu}}\left(\frac{x}{y}\right)^{\nu} = \frac{y^2}{1-x^2}\left(\frac{y}{x}\right)^\nu.
\end{split}\end{equation}

For $\nu = 0,1,2$, it is straightforward to show that $LHS(\nu) = RHS(\nu)$; %indeed, 
%\begin{equation}\begin{split}\notag
%LHS(0) &=  A_0 +\left(1+\frac{x^2}{y^2}\right)A_{0,0} -\frac{x}{y}A_{0,0}U_1
%= 1-y^2 + \left(1+\frac{x^2}{y^2}\right)\frac{y^2}{1-x^2}-\frac{x}{y} \frac{y^2}{1-x^2} U_1 \\
%& = 1-y^2 + \frac{y^2}{1-x^2} + \frac{x^2}{1-x^2} - \frac{xy}{1-x^2} U_1
%= \frac{1-x^2y^2}{1-x^2} - \frac{xy}{1-x^2} U_1 
%= RHS(0)
%\end{split}\end{equation}
%and similarly (trivial computations)
%\begin{equation}\begin{split}\notag
%LHS(1) &=  \tfrac{x}{y} \{ A_{1} +(1+\tfrac{x^2}{y^2}) ( A_{0,1}  +  A_{1,1} U_1) -\tfrac{x}{y}( A_{0,1}U_1 + A_{1,1}U_1^2 )\} = RHS(1)\\
%LHS(2) &=  (\tfrac{x}{y})^2 \{ A_{2} +(1+\tfrac{x^2}{y^2}) ( A_{0,2}  +  A_{1,2} U_1 +  A_{2,2}U_2)  -\tfrac{x}{y}( A_{0,2}U_1 + A_{1,2}U_1^2 + A_{2,2}U_1U_2 )\} = RHS(2).
%\end{split}\end{equation}
hence, from now on we can assume that $\nu\geq 3$. 
%Since
%\begin{equation}\begin{split}\notag
%& U_1U_i = U_{i+1} + U_{i-1} \quad \text{ for }i>0, \\
%& \sum_{i=0}^\nu A_{i,\nu}U_1 U_i = A_{0,\nu}U_1 + \sum_{i=1}^\nu A_{i,\nu}(U_{i-1}+U_{i+1})  , \\
%& \sum_{i=1}^{\nu} A_{i,\nu}U_{i-1} = \sum_{i=0}^{\nu-1}A_{i+1,\nu}U_i 
%= A_{1,\nu} + A_{2,\nu} U_1 + \sum_{i=2}^{\nu-2}A_{i+1,\nu}U_i + A_{\nu,\nu}U_{\nu-1},\\
%& \sum_{i=1}^{\nu} A_{i,\nu}U_{i+1} = \sum_{i=2}^{\nu+1} A_{i-1.\nu}U_i 
%= \sum_{i=2}^{\nu-2} A_{i-1,\nu}U_i + A_{\nu-2,\nu}U_{\nu-1} + A_{\nu-1,\nu}U_{\nu} + A_{\nu,\nu}U_{\nu+1} ,
%\end{split}\end{equation}
Using \eqref{recurrence}, it is easy to get
\begin{equation}\begin{split}\notag
LHS(\nu)
%&=\left(\frac{x}{y}\right)^\nu \bigg\{ A_{\nu} +\left(1+\frac{x^2}{y^2}\right) \bigg[A_{0,\nu} + A_{1,\nu}U_1 + \sum_{i=2}^{\nu -2}A_{i,\nu}U_i + A_{\nu-1,\nu}U_{\nu-1} + A_{\nu,\nu} U_\nu\bigg] \\
%&\quad -\frac{x}{y}A_{0,\nu}U_1  -\frac{x}{y}\bigg[  A_{1,\nu} + A_{2,\nu} U_1 + \sum_{i=2}^{\nu-2}A_{i+1,\nu}U_i +  A_{\nu,\nu}U_{\nu-1} \bigg]\\
%&\quad-\frac{x}{y}\bigg[ \sum_{i=2}^{\nu-2} A_{i-1,\nu}U_i + A_{\nu-2,\nu}U_{\nu-1}  + A_{\nu-1,\nu}U_{\nu} + A_{\nu,\nu}U_{\nu+1} \bigg] \bigg\} \\
&=\left(\frac{x}{y}\right)^\nu \bigg\{ A_\nu +\left(1+\frac{x^2}{y^2}\right) \bigg[A_{0,\nu} + A_{1,\nu}U_1 + A_{\nu-1,\nu}U_{\nu-1}  + A_{\nu,\nu} U_\nu\bigg]\\
&\quad -\frac{x}{y}\bigg[ A_{0\nu}U_1 + A_{1,\nu} + A_{2,\nu} U_1 + A_{\nu,\nu}U_{\nu-1} + A_{\nu-2,\nu}U_{\nu-1}  + A_{\nu-1,\nu}U_{\nu} + A_{\nu,\nu}U_{\nu+1} \bigg] \\
&\quad+\sum_{i=2}^{\nu -2} \bigg[ \left(1+\frac{x^2}{y^2}\right)A_{i,\nu} -\frac{x}{y}A_{i-1,\nu}   -\frac{x}{y} A_{i-1,\nu}\bigg] U_i  \bigg\}.
\end{split}\end{equation}

\begin{lemma}\label{telescopiclemma}
For any positive integers $\nu,i$ with $0\leq i< \nu-1$,
%\todo{maybe remove this condition?}
we have
$$\left(1+\frac{x^2}{y^2}\right)A_{i,\nu} -\frac{x}{y}A_{i-1,\nu}-\frac{x}{y} A_{i-1,\nu}=0. $$
\end{lemma}
\proof 
By definition of $A_{i,\nu}$, denoting $z=y/x$, the left hand side (say $f$) reads
\begin{equation}\begin{split}\notag
f(\nu,i) = &
(1+z^{-2})y^2\bigg( -z^i -(y^{-2}-z^{-2})z^{-i}\sum_{i+1\leq c\leq \nu}z^{2c} + z^{2\nu-i}  \bigg)\\
&\quad -z^{-1}y^2\bigg(  -z^{i+1} -(y^{-2}-z^{-2})z^{-i-1}\sum_{i+2\leq c\leq \nu}z^{2c} + z^{2\nu-i-1}  \bigg)\\
&\quad -z^{-1}y^2\bigg(  -z^{i-1} -(y^{-2}-z^{-2})z^{-i+1}\sum_{i\leq c\leq \nu}z^{2c} + z^{2\nu-i+1}  \bigg).
\end{split}\end{equation}
We now split the factor $y^{-2}-z^{-2}$; first we consider the contribution coming from $z^{-2}$ (which we will denote $f^\flat$), the other ($f^\sharp$) will be handled below. We have
\begin{equation}\begin{split}\notag
f^\flat(\nu,i) &
%(1+z^{-2})y^2\bigg( -z^i +z^{-i-2}\sum_{i+1\leq c\leq \nu}z^{2c} + z^{2\nu-i}  \bigg)\\
%&\quad -z^{-1}y^2\bigg(  -z^{i+1} + z^{-i-3}\sum_{i+2\leq c\leq \nu}z^{2c} + z^{2\nu-i-1}  \bigg)\\
%&\quad -z^{-1}y^2\bigg(  -z^{i-1} + z^{-i-1}\sum_{i\leq c\leq \nu}z^{2c} + z^{2\nu-i+1}  \bigg)\\
%&=(1+z^{-2})y^2\bigg( -z^i +z^{-i-2}\bigg( z^{2i+2}+z^{2i+4}+\dots+z^{2\nu-2}+z^{2\nu} \bigg) + z^{2\nu-i} \bigg)\\
%&\quad -y^2\bigg(  -z^{i} + z^{-i-4}\bigg( z^{2i+4}+z^{2i+6}+\dots+z^{2\nu-2}+z^{2\nu} \bigg)  + z^{2\nu-i-2}  \bigg)\\
%&\quad -z^{-2}y^2\bigg(  -z^{i} + z^{-i}\bigg( z^{2i}+z^{2i+2}+\dots+z^{2\nu-2}+z^{2\nu} \bigg)  + z^{2\nu-i+2}  \bigg)\\
%&=(1+z^{-2})y^2\bigg( -z^i + z^{i}+z^{i+2}+\dots+z^{2\nu-i-4}+z^{2\nu-i-2} + z^{2\nu-i} \bigg)\\
%&\quad -y^2\bigg(  -z^{i} +  z^{i}+z^{i+2}+\dots+z^{2\nu-i-6}+z^{2\nu-i-4} + z^{2\nu-i-2}  \bigg)\\
%&\quad -z^{-2}y^2\bigg(  -z^{i} +  z^{i}+z^{i+2}+\dots+z^{2\nu-i-2}+z^{2\nu-i}   + z^{2\nu-i+2}  \bigg)\\
=(1+z^{-2})y^2( z^{i+2}+\dots+z^{2\nu-i-2} + z^{2\nu-i} )\\
&\quad -y^2(  z^{i+2}+\dots+ z^{2\nu-i-2}  )\\
&\quad -z^{-2}y^2(  z^{i+2}+\dots+z^{2\nu-i-2}+z^{2\nu-i}   + z^{2\nu-i+2}  )\\
&=(1+z^{-2})y^2( z^{2\nu-i} )-z^{-2}y^2(  z^{2\nu-i}   + z^{2\nu-i+2} ) 
= 0.
\end{split}\end{equation}
We finally treat $f^\sharp$ very similarly, as
\begin{equation}\begin{split}\notag
f^\sharp(\nu,i) %&= 
%(1+z^{-2})y^2\bigg( -z^i -y^{-2}z^{-i}\bigg( z^{2i+2}+z^{2i+4}+\dots+z^{2\nu-2}+z^{2\nu} \bigg) + z^{2\nu-i}  \bigg)\\
%&\quad -y^2\bigg(  -z^{i} -y^{-2} z^{-i-2}\bigg( z^{2i+4}+z^{2i+6}+\dots+z^{2\nu-2}+z^{2\nu} \bigg)  + z^{2\nu-i-2}  \bigg)\\
%&\quad -z^{-2}y^2\bigg(  -z^{i} -y^{-2} z^{-i+2}\bigg( z^{2i}+z^{2i+2}+\dots+z^{2\nu-2}+z^{2\nu} \bigg) + z^{2\nu-i+2}  \bigg)\\
%&=(1+z^{-2})y^2\bigg( -z^i -y^{-2}\bigg( z^{i+2}+\dots+z^{2\nu-i-2}+z^{2\nu-i} \bigg) + z^{2\nu-i}  \bigg)\\
%&\quad -y^2\bigg(  -z^{i} -y^{-2} \bigg( z^{i+2}+\dots+z^{2\nu-i-2} \bigg)  + z^{2\nu-i-2}  \bigg)\\
%&\quad -z^{-2}y^2\bigg(  -z^{i} -y^{-2} \bigg( z^{i+2}+\dots+z^{2\nu-i-2}+z^{2\nu-i}+z^{2\nu-i+2} \bigg) + z^{2\nu-i+2}  \bigg)\\
&=(1+z^{-2})y^2(  -y^{-2} z^{2\nu-i}  + z^{2\nu-i}  )\\
&\quad -y^2(   z^{2\nu-i-2}  )\\
&\quad -z^{-2}y^2(   -y^{-2} ( z^{2\nu-i}+z^{2\nu-i+2} ) + z^{2\nu-i+2}  ) 
= 0.
\end{split}\end{equation}
\endproof

Thanks to Lemma \ref{telescopiclemma}, we obtain
\begin{equation}\begin{split}\label{6may.1}\notag
LHS(\nu)
&=\left(\frac{x}{y}\right)^\nu \bigg\{ A_\nu +\left(1+\frac{x^2}{y^2}\right) \bigg[A_{0,\nu} + A_{1,\nu}U_1 + A_{\nu-1,\nu}U_{\nu-1}+ A_{\nu,\nu} U_\nu\bigg]\\
&\hspace{0.1cm}
 -\frac{x}{y}\bigg[ A_{0,\nu}U_1 + A_{1,\nu} + A_{2,\nu} U_1 + A_{\nu,\nu}U_{\nu-1} + A_{\nu-2,\nu}U_{\nu-1} + A_{\nu-1,\nu}U_{\nu} + A_{\nu,\nu}U_{\nu+1} \bigg]  \bigg\}\\
&=: B_{0} + B_1U_1 + B_{\nu-1}U_{\nu-1} + B_{\nu}U_{\nu} + B_{\nu+1}U_{\nu+1}.
\end{split}\end{equation}

To show that $LHS(\nu) = RHS(\nu)$ for all $\nu\geq 3$, it remains to show that the coefficients (denoted by $B_i$) of $U_0=1, U_1, U_{\nu-1}, U_{\nu}, U_{\nu+1}$ of $LHS$ equal those of $RHS$. From 
the above display %\eqref{6may.1}, 
we see that the coefficient of $U_{\nu +1}$ in $LHS$ is
$$ B_{\nu+1} = -\left(\frac{x}{y}\right)^\nu \frac{x}{y} A_{\nu,\nu} 
= -\left(\frac{x}{y}\right)^\nu \frac{x}{y} \frac{y^2}{1-x^2}\left(\frac{y}{x}\right)^\nu 
=   -  \frac{xy}{1-x^2}  $$
and then equals the coefficient of $U_{\nu+1}$ of $RHS$. Similarly, the coefficient of $U_\nu$ of $LHS$ is
$$ B_\nu=\left(\frac{x}{y}\right)^\nu\left( \left(1+\frac{x^2}{y^2}\right) A_{\nu,\nu}
 -\frac{x}{y} A_{\nu-1,\nu}\right)  
= \frac{1+x^2y^2}{1-x^2} $$
%\begin{equation}\begin{split}\notag 
%B_\nu=&\left(\frac{x}{y}\right)^\nu\left( \left(1+\frac{x^2}{y^2}\right) A_{\nu,\nu}
% -\frac{x}{y} A_{\nu-1,\nu}\right) \\ 
% =& (\tfrac{x}{y})^\nu \left( (1+\tfrac{x^2}{y^2}) \tfrac{y^2}{1-x^2}\left(\tfrac{y}{x}\right)^{\nu}
%-\tfrac{x}{y} \left(
%-y^2\left(\tfrac{y}{x}\right)^{\nu-1} - (1-x^2)(\tfrac{x}{y})^{\nu-1}\left(\tfrac{y}{x}\right)^{2\nu} + \tfrac{y^{2(\nu+1)}}{x^{2\nu}}(\tfrac{x}{y})^{\nu-1}\right)\right)\\
%=& \frac{1+x^2y^2}{1-x^2}
%\end{split}\end{equation}
as in $RHS$ and
\begin{equation}\begin{split}\notag 
B_{\nu-1}=&\left(\frac{x}{y}\right)^\nu\left( \left(1+\frac{x^2}{y^2}\right)A_{\nu-1,\nu}-\frac{x}{y} A_{\nu,\nu}
 -\frac{x}{y}  A_{\nu-2,\nu}\right) = -\frac{xy}{1-x^2}.
\end{split}\end{equation}
Furthermore, for $U_1$ the coefficient in $LHS$ is 
\begin{equation}\label{B1}\notag
B_1 = \left(\frac{x}{y}\right)^\nu \bigg( \left(1+\frac{x^2}{y^2}\right)  A_{1,\nu}
 -\frac{x}{y} A_{0,\nu} -\frac{x}{y}A_{2,\nu} \bigg) = 0
 \end{equation}
 %\todo[inline]{da dimostrare}
 like in $RHS$ (it follows from Lemma \ref{telescopiclemma} as a special case). Finally, with the same strategy as in Lemma \ref{telescopiclemma}, we can show that the constant term is
\begin{equation}\label{B0}\notag
 B_0 = \left(\frac{x}{y}\right)^\nu \bigg\{ A_\nu +\left(1+\frac{x^2}{y^2}\right) A_{0,\nu}  -\frac{x}{y} A_{1,\nu} \bigg\}=0.
\end{equation}
%, as shown below. 
%\begin{lemma}\label{lemma2}
%In the notations above, we have
%$$A_\nu +\left(1+\frac{x^2}{y^2}\right) A_{0,\nu}  -\frac{x}{y} A_{1,\nu}=0.$$
%\end{lemma}
%\proof
%Denoting $z=y/x$ and $\mathcal S=\sum_{c=1}^\nu z^{2c}$, we can write the left hand side as
%\begin{equation}\begin{split}\notag
%%A_\nu +(1+\tfrac{x^2}{y^2}) A_{0,\nu}  -\tfrac{x}{y} A_{1,\nu} 
%&= 1+(1-x^2)\mathcal S - y^2z^{2\nu}\\
%& \hspace{1cm}+(1+z^{-2})( -y^2-(1-x^2)\mathcal S + y^2 z^{2\nu} )\\
%& \hspace{1cm}-z^{-1}( -y^2z-(1-x^2)z^{-1}(\mathcal S-z^2)+ y^2 z^{2\nu}z^{-1} )\\
%&= 1+(1-x^2)\mathcal S - y^2z^{2\nu}\\
%& \hspace{1cm} -y^2-(1-x^2)\mathcal S + y^2 z^{2\nu}  + z^{-2}( -y^2-(1-x^2)\mathcal S + y^2 z^{2\nu} )\\
%& \hspace{1cm} + y^2 + z^{-2}(1-x^2)\mathcal S - (1-x^2) - y^2 z^{2\nu}z^{-2}\\
%&= 1+(1-x^2)\mathcal S - y^2z^{2\nu}\\
%& \hspace{1cm} -y^2-(1-x^2)\mathcal S + y^2 z^{2\nu}   
%-z^{-2}y^2-z^{-2}(1-x^2)\mathcal S + y^2 z^{2\nu-2} \\
%& \hspace{1cm} + y^2 + z^{-2}(1-x^2)\mathcal S - (1-x^2) - y^2 z^{2\nu-2}\\
%&= 1 - y^2z^{2\nu} -y^2 + y^2 z^{2\nu}   
%- x^2 + y^2 z^{2\nu-2} + y^2  - 1 + x^2 - y^2 z^{2\nu-2}\\
%& =0.
%\end{split}\end{equation}
\endproof

\medskip
\textbf{Acknowledgments}. We would like to thank Sieg Baluyot for many helpful discussions and the referee for spotting a critical mistake in an earlier version of the paper. This work is supported by the FRG grant DMS 1854398. The second author is member of the INdAM group GNAMPA.

{\small

\end{document}